\documentclass[a4paper]{amsart}
\usepackage{amssymb,enumerate,mathrsfs,mathtools,stmaryrd,nicefrac}
\usepackage[all]{xy}

\title{Categorical aspects of bivariant K-theory}

\author{Ralf Meyer}

\address{Mathematisches Institut\\
  Georg-August-Universit\"at G\"ottingen\\
  Bunsenstra\ss{}e 3--5\\
  37073 G\"ottingen\\
  Deutschland}
\email{rameyer at uni-math.gwdg.de}

\usepackage[pdftitle={Categorical aspects of bivariant K-theory},
  pdfauthor={Ralf Meyer},
  pdfsubject={Mathematics; MSC 19K35, 46L80}
]{hyperref}
\usepackage[lite]{amsrefs}

\usepackage{microtype}

\theoremstyle{plain}
\newtheorem{theorem}{Theorem}
\newtheorem{lemma}[theorem]{Lemma}
\newtheorem{proposition}[theorem]{Proposition}
\newtheorem{corollary}[theorem]{Corollary}
\theoremstyle{definition}
\newtheorem{definition}[theorem]{Definition}
\theoremstyle{remark}
\newtheorem{remark}[theorem]{Remark}
\newtheorem{example}[theorem]{Example}

\newcommand*{\Cstar}{\texorpdfstring{$C^*$\nobreakdash-\hspace{0pt}}{C*-}}
\newcommand*{\Star}{\texorpdfstring{$^*$\nobreakdash-\hspace{0pt}}{*-}}

\newcommand*{\C}{\mathbb C}
\newcommand*{\Z}{\mathbb Z}
\newcommand*{\Ztwo}{\mathbb Z/2}
\newcommand*{\R}{\mathbb R}
\newcommand*{\N}{\mathbb N}

\newcommand*{\Mat}{\mathbb{M}}
\newcommand*{\Comp}{{\mathbb{K}}}
\newcommand*{\Bound}{{\mathbb{B}}}
\newcommand*{\Sphere}{\mathbb S}

\newcommand*{\Cat}{\mathfrak{C}}

\newcommand*{\Univ}{\mathfrak{Univ}}

\newcommand*{\Cstarcat}{\mathfrak{C^*alg}}
\newcommand*{\CstarcatG}[1][G]{#1\text{-}\mathfrak{C^*alg}}
\newcommand*{\Cstarsep}{\mathfrak{C^*sep}}
\newcommand*{\CstarsepG}[1][G]{#1\text{-}\mathfrak{C^*sep}}
\newcommand*{\Cstarnice}{\mathfrak{S}}
\newcommand*{\Corr}[1]{\mathfrak{Corr}_{#1}}
\newcommand*{\KKcat}{\mathfrak{KK}}
\newcommand*{\Ecat}{\mathfrak{E}}

\newcommand*{\Cfull}{C^*}
\newcommand*{\Cred}{C^*_{\mathrm{red}}}

\newcommand*{\KK}{\mathrm{KK}}
\newcommand*{\RKK}{\mathrm{RKK}}
\newcommand*{\K}{\mathrm{K}}
\newcommand*{\E}{\mathrm{E}}

\newcommand*{\ev}{\mathrm{ev}}
\newcommand*{\const}{\mathrm{const}}
\newcommand*{\Ho}{\mathrm{Ho}}

\newcommand*{\Cyl}{\mathrm{Cyl}}
\newcommand*{\Sus}{\mathrm{Sus}}
\newcommand*{\Cone}{\mathrm{Cone}}

\newcommand*{\ID}{{\mathrm{id}}}

\DeclareMathOperator{\Ind}{Ind}
\DeclareMathOperator{\Res}{Res}

\DeclareMathOperator{\Hom}{Hom}
\DeclareMathOperator{\Ext}{Ext}
\DeclareMathOperator{\Map}{Map}
\DeclareMathOperator{\Asymp}{Asymp}

\newcommand*{\CONT}{{\mathcal{C}}}
\newcommand*{\Hils}{{\mathcal{H}}}
\newcommand*{\Hilm}[1][E]{{\mathcal{#1}}}
\newcommand*{\EG}{{\mathcal{E}}}
\newcommand*{\Fam}{{\mathcal{F}}}
\newcommand*{\genci}{\langle\mathcal{CI}\rangle}
\newcommand*{\CC}{\mathcal{CC}}

\newcommand*{\Mult}{{\mathcal{M}}}

\newcommand*{\brd}{-\hspace{0pt}}
\newcommand*{\nbd}{\nobreakdash-\hspace{0pt}}

\newcommand*{\norm}[1]{\lVert#1\rVert}
\newcommand*{\inOb}{\mathrel{\in\!\in}}         
\newcommand*{\barotimes}{\mathbin{\bar{\otimes}}}
\newcommand*{\minotimes}{\mathbin{\otimes_{\mathrm{min}}}}
\newcommand*{\maxotimes}{\mathbin{\otimes_{\mathrm{max}}}}

\newcommand*{\rcross}{\mathbin{\ltimes_{\mathrm{r}}}}
\newcommand*{\cross}{\mathbin{\ltimes}}
\newcommand*{\defeq}{\mathrel{\vcentcolon=}}
\newcommand*{\into}{\rightarrowtail}
\newcommand*{\prto}{\twoheadrightarrow}
\newcommand*{\blank}{\text{\textvisiblespace}}

\newcommand*{\pt}{\star}

\begin{document}

\begin{abstract}
  This survey article on bivariant \textsc{Kasparov} theory and
  \(\E\)\nbd{}theory is mainly intended for readers with a background in
  homotopical algebra and category theory.  We approach both bivariant
  \(\K\)\nbd{}theories via their universal properties and equip them
  with extra structure such as a tensor product and a triangulated
  category structure.  We discuss the construction of the
  \textsc{Baum}--\textsc{Connes} assembly map via localisation of
  categories and explain how this is related to the purely topological
  construction by \textsc{Davis} and \textsc{L\"uck}.
\end{abstract}




\maketitle

\section{Introduction}
\label{sec:intro}

Non-commutative topology deals with topological properties of
\Cstar{}algebras.  Already in the 1970s, the classification of
AF-algebras by \(\K\)\nbd{}theoretic data \cite{Elliott:Classify-AF} and
the work of \textsc{Brown}--\textsc{Douglas}--\textsc{Fillmore} on
essentially normal operators \cite{BDF} showed clearly that topology
provides useful tools to study \Cstar{}algebras.  A breakthrough was
\textsc{Kasparov}'s construction of a bivariant \(\K\)\nbd{}theory for
separable \Cstar{}algebras.  Besides its applications within
\Cstar{}algebra theory, it also yields results in classical topology
that are hard or even impossible to prove without it.  A typical example
is the \textsc{Novikov} conjecture, which deals with the homotopy
invariance of certain invariants of smooth manifolds with a given
fundamental group.  This conjecture has been verified for many groups
using \textsc{Kasparov} theory, starting with \cite{Kasparov:Novikov}.
The \Cstar{}algebraic formulation of the \textsc{Novikov} conjecture is
closely related to the \textsc{Baum}--\textsc{Connes} conjecture, which
deals with the computation of the \(\K\)\nbd{}theory \(\K_*(\Cred G)\)
of reduced group \Cstar{}algebras and has been one of the centres of
attention in non-commutative topology in recent years.

The \textsc{Baum}--\textsc{Connes} conjecture in its original
formulation \cite{Baum-Connes-Higson} only deals with a single
\(\K\)\nbd{}theory group; but a better understanding requires a
different point of view.  The approach by \textsc{Davis} and
\textsc{L\"uck} in \cite{Davis-Lueck:Assembly} views it as a natural
transformation between two homology theories for
\(G\)\nbd{}CW-complexes.  An analogous approach in the \Cstar{}algebra
framework appeared in \cite{Meyer-Nest:BC}.  These approaches to the
\textsc{Baum}--\textsc{Connes} conjecture show the importance of
studying not just single \Cstar{}algebras, but categories of
\Cstar{}algebras and their properties.  Older ideas like the universal
property of \textsc{Kasparov} theory are of the same nature.  Studying
categories of objects instead of individual objects is becoming more and
more important in algebraic topology and algebraic geometry as well.

Several mathematicians have suggested, therefore, to apply general
constructions with categories (with additional structure) like
generators, \textsc{Witt} groups, the centre, and support varieties to
the \Cstar{}algebra context.  Despite the warning below, this seems a
promising project, where little has been done so far.  To prepare for
this enquiry, we summarise some of the known properties of categories of
\Cstar{}algebras; we cover tensor products, some homotopy theory,
universal properties, and triangulated structures.  In addition, we
examine the Universal Coefficient Theorem and the
\textsc{Baum}--\textsc{Connes} assembly map.

Despite many formal similarities, the homotopy theory of spaces and
non\brd{}commutative topology have a very different focus.

On the one hand, most of the complexities of the stable homotopy
category of spaces vanish for \Cstar{}algebras because only very few
homology theories for spaces have a non-commutative counterpart: any
functor on \Cstar{}algebras satisfying some reasonable assumptions must
be closely related to \(\K\)\nbd{}theory.  Thus special features of
topological \(\K\)\nbd{}theory become more transparent when we work with
\Cstar{}algebras.

On the other hand, analysis may create new difficulties, which appear to
be very hard to study topologically.  For instance, there exist
\Cstar{}algebras with vanishing \(\K\)\nbd{}theory which are
nevertheless non-trivial in \textsc{Kasparov} theory; this means that
the Universal Coefficient Theorem fails for them.  I know no non-trivial
topological statement about the subcategory of the \textsc{Kasparov}
category consisting of \Cstar{}algebras with vanishing
\(\K\)\nbd{}theory; for instance, I know no compact objects.

It may be necessary, therefore, to restrict attention to suitable
``bootstrap'' categories in order to exclude pathologies that have
nothing to do with classical topology.  More or less by design, the
resulting categories will be localisations of purely topological
categories, which we can also construct without mentioning
\Cstar{}algebras.  For instance, we know that the
\textsc{Rosenberg}--\textsc{Schochet} bootstrap category is equivalent
to a full subcategory of the category of \(BU\)-module spectra.  But we
can hope for more interesting categories when we work equivariantly with
respect to, say, discrete groups.

\section{Additional structure in \texorpdfstring{$C^*$-}{C*-}algebra
  categories}
\label{sec:Cstar_structure}

We assume that the reader is familiar with some basic properties of
\Cstar{}algebras, including the definition (see for instance
\cites{Arveson:Invitation, Davidson:Cstar_example}).  As usual, we allow
non-unital \Cstar{}algebras.  We define some categories of
\Cstar{}algebras in~\textsection\ref{sec:categories_Cstar} and consider
group \Cstar{}algebras and crossed products
in~\textsection\ref{sec:group_crossed}.  Then we discuss \Cstar{}tensor
products and mention the notions of nuclearity and exactness in
\textsection\ref{sec:tensor}.  The upshot is that \(\Cstarcat\) and
\(\CstarcatG\) carry two structures of symmetric monoidal category,
which coincide for nuclear \Cstar{}algebras.  We prove in
\textsection\ref{sec:limits} that \(\Cstarcat\) and \(\CstarcatG\) are
bicomplete, that is, all diagrams in them have both a limit and a
colimit.  We equip morphism spaces between \Cstar{}algebras with a
canonical base point and topology in \textsection\ref{sec:enrichment};
thus the category of \Cstar{}algebras is enriched over the category of
pointed topological spaces.  In \textsection\ref{sec:cylinders}, we
define mapping cones and cylinders in categories of \Cstar{}algebras;
these rudimentary tools suffice to carry over some basic homotopy
theory.

\subsection{Categories of \texorpdfstring{$C^*$-}{C*-}algebras}
\label{sec:categories_Cstar}

\begin{definition}
  \label{def:Cstarcat}
  The category of \Cstar{}algebras is the category \(\Cstarcat\) whose
  objects are the \Cstar{}algebras and whose morphisms \(A\to B\) are
  the \Star{}homomorphisms \(A\to B\); we denote this set of morphisms
  by \(\Hom(A,B)\).

  A \Cstar{}algebra is called \emph{separable} if it has a countable
  dense subset.  We often restrict attention to the full subcategory
  \(\Cstarsep\subseteq\Cstarcat\) of separable \Cstar{}algebras.
\end{definition}

Examples of \Cstar{}algebras are group \Cstar{}algebras and
\Cstar{}crossed products.  We briefly recall some relevant properties of
these constructions.  A more detailed discussion can be found in many
textbooks such as~\cite{Pedersen:Automorphism}.

\begin{definition}
  \label{def:inOb}
  We write \(A\inOb\Cat\) to denote that~\(A\) is an object of the
  category~\(\Cat\).  The notation \(f\in\Cat\) means that~\(f\) is a
  morphism in~\(\Cat\); but to avoid confusion we always specify domain
  and target and write \(f\in\Cat(A,B)\) instead of \(f\in\Cat\).
\end{definition}

\subsection{Group actions, and crossed products}
\label{sec:group_crossed}

For any locally compact group~\(G\), we have a \emph{reduced group
  \Cstar{}algebra} \(\Cred(G)\) and a \emph{full group
  \Cstar{}algebra} \(\Cfull(G)\).  Both are defined as completions of
the group \textsc{Banach} algebra \((L^1(G),*)\) for suitable
\Cstar{}norms and are related by a canonical surjective
\Star{}homomorphism \(\Cfull(G)\to\Cred(G)\), which is an isomorphism
if and only if~\(G\) is \emph{amenable}.

The norm on \(\Cfull(G)\) is the maximal \Cstar{}norm, so that any
strongly continuous unitary representation of~\(G\) on a
\textsc{Hilbert} space induces a \Star{}representation of \(\Cfull(G)\).
The norm on \(\Cred(G)\) is defined using the regular representation
of~\(G\) on \(L^2(G)\); hence a representation of~\(G\) only induces a
\Star{}representation of \(\Cred(G)\) if it is weakly contained in the
regular representation.  For reductive \textsc{Lie} groups and reductive
\(p\)\nbd{}adic groups, these representations are exactly the tempered
representations, which are much easier to classify than all unitary
representations.

\begin{definition}
  \label{def:Cstarcat_G}
  A \emph{\(G\)\nbd{}\Cstar{}algebra} is a \Cstar{}algebra~\(A\) with a
  strongly continuous representation of~\(G\) by \Cstar{}algebra
  automorphisms.  The category of \(G\)\nbd{}\Cstar{}algebras is the
  category \(\CstarcatG\) whose objects are the
  \(G\)\nbd{}\Cstar{}algebras and whose morphisms \(A\to B\) are the
  \(G\)\nbd{}equivariant \Star{}homomorphisms \(A\to B\); we denote this
  morphism set by \(\Hom_G(A,B)\).
\end{definition}

\begin{example}
  \label{exa:Z_action}
  If \(G=\Z\), then a \(G\)-\Cstar{}algebra is nothing but a pair
  \((A,\alpha)\) consisting of a \Cstar{}algebra~\(A\) and a
  \Star{}automorphism \(\alpha\colon A\to A\): let~\(\alpha\) be the
  action of the generator \(1\in\Z\).
\end{example}

Equipping \Cstar{}algebras with a trivial action provides a functor
\begin{equation}
  \label{eq:trivial_action}
  \tau\colon \Cstarcat\to\CstarcatG, \qquad A\mapsto A_\tau.
\end{equation}
Since~\(\C\) has only the identity automorphism, the trivial action is
the only way to turn~\(\C\) into a \(G\)-\Cstar{}algebra.

The \emph{full and reduced \Cstar{}crossed products} are versions of
the full and reduced group \Cstar{}algebras with coefficients in
\(G\)-\Cstar{}algebras (see~\cite{Pedersen:Automorphism}).  They
define functors
\[
G\cross\blank,G\rcross\blank\colon \CstarcatG\to\Cstarcat,
\qquad
A\mapsto G\cross A,\quad
G\rcross A,
\]
such that \(G\cross\C=\Cfull(G)\) and \(G\rcross\C=\Cred(G)\).

\begin{definition}
  \label{def:extension}
  A diagram \(I\to E\to Q\) in \(\Cstarcat\) is an \emph{extension} if
  it is isomorphic to the canonical diagram \(I\to A \to A/I\) for some
  ideal~\(I\) in a \Cstar{}algebra~\(A\); extensions in \(\CstarcatG\)
  are defined similarly, using \(G\)\nbd{}invariant ideals in
  \(G\)-\Cstar{}algebras.  We write \(I\into E\prto Q\) to denote
  extensions.
\end{definition}

Although \Cstar{}algebra extensions have some things in common with
extensions of, say, modules, there are significant differences because
\(\Cstarcat\) is not \textsc{Abel}ian, not even additive.

\begin{proposition}
  \label{pro:fullcross_exact}
  The full crossed product functor \(G\cross\blank\colon
  \CstarcatG\to\Cstarcat\) is \emph{exact} in the sense that it maps
  extensions in \(\CstarcatG\) to extensions in \(\Cstarcat\).
\end{proposition}

\begin{proof}
  This is Lemma 4.10 in \cite{Guentner-Higson-Trout}.
\end{proof}

\begin{definition}
  \label{def:exact_group}
  A locally compact group~\(G\) is called \emph{exact} if the reduced
  crossed product functor \(G\rcross\blank\colon
  \CstarcatG\to\Cstarcat\) is exact.
\end{definition}

Although this is not apparent from the above definition,
\emph{exactness} is a geometric property of a group: it is equivalent
to \textsc{Yu}'s property~(A) or to the existence of an amenable
action on a compact space \cite{Ozawa:Amenable_action}.

Most groups you know are exact.  The only source of non-exact groups
known at the moment are \textsc{Gromov}'s random groups.  Although
exactness might remind you of the notion of flatness in homological
algebra, it has a very different flavour.  The difference is that the
functor \(G\rcross\blank\) always preserves injections and surjections.
What may go wrong for non-exact groups is exactness in the \emph{middle}
(compare the discussion before
Proposition~\ref{pro:exact_on_cpsplit_tensor}).  Hence we cannot study
the lack of exactness by derived functors.

Even for non-exact groups, there is a class of extensions for which
reduced crossed products are always exact:

\begin{definition}
  \label{def:cp-split}
  A \emph{section} for an extension
  \begin{equation}
    \label{eq:extension_cp-split}
    I\overset{i}\into E \overset{p}\prto Q
  \end{equation}
  in~\(\CstarcatG\) is a map (of sets) \(s\colon Q\to E\) with \(p\circ
  s=\ID_Q\).  We call~\eqref{eq:extension_cp-split} \emph{split} if
  there is a section that is a \(G\)\nbd{}equivariant
  \Star{}homomorphism.  We call~\eqref{eq:extension_cp-split}
  \emph{\(G\)\nbd{}equivariantly cp-split} if there is a
  \(G\)\nbd{}equivariant, completely positive, contractive, linear
  section.
\end{definition}

Sections are also often called \emph{lifts}, \emph{liftings}, or
\emph{splittings}.

\begin{proposition}
  \label{pro:exact_on_cpsplit}
  Both the reduced and the full crossed product functors map split
  extensions in \(\CstarcatG\) again to split extensions
  in~\(\Cstarcat\) and \(G\)\nbd{}equivariantly cp-split extensions in
  \(\CstarcatG\) to cp-split extensions in~\(\Cstarcat\).
\end{proposition}

\begin{proof}
  Let \(K\overset{i}\into E\overset{p}\prto Q\) be an extension in
  \(\CstarcatG\).  Proposition~\ref{pro:fullcross_exact} shows that
  \(G\cross K\into G\cross E\prto G\cross Q\) is again an extension.
  Since reduced and full crossed products are functorial for equivariant
  completely positive contractions, this extension is split or cp-split
  if the original extension is split or equivariantly cp-split,
  respectively.  This yields the assertions for full crossed products.

  Since a \Star{}homomorphism with dense range is automatically
  surjective, the induced map \(G\rcross p\colon G\rcross E\to G\rcross
  Q\) is surjective.  It is evident from the definition of reduced
  crossed products that \(G\rcross i\) is injective.  What is unclear is
  whether the range of \(G\rcross i\) and the kernel of \(G\rcross p\)
  coincide.  As for the full crossed product, a \(G\)\nbd{}equivariant
  completely positive contractive section \(s\colon Q\to E\) induces a
  completely positive contractive section \(G\rcross s\) for \(G\rcross
  p\).  The linear map
  \[
  \varphi\defeq \ID_{G\rcross E} - (G\rcross s)\circ (G\rcross p)\colon
  G\rcross E\to G\rcross E
  \]
  is a retraction from \(G\rcross E\) onto the kernel of \(G\rcross p\)
  by construction.  Furthermore, it maps the dense subspace \(L^1(G,E)\)
  into \(L^1(G,K)\).  Hence it maps all of \(G\rcross E\) into
  \(G\rcross K\).  This implies \(G\rcross K= \ker(G\rcross p)\) as
  desired.
\end{proof}

\subsection{Tensor products and nuclearity}
\label{sec:tensor}

Most results in this section are proved in detail in
\cites{Murphy:book, Wassermann:Exact}.  Let \(A_1\) and~\(A_2\) be two
\Cstar{}algebras.  Their (algebraic) tensor product \(A_1\otimes A_2\)
is still a \Star{}algebra.  A \emph{\Cstar{}tensor product} of \(A_1\)
and~\(A_2\) is a \Cstar{}completion of \(A_1\otimes A_2\), that is, a
\Cstar{}algebra that contains \(A_1\otimes A_2\) as a dense
\Star{}subalgebra.  A \Cstar{}tensor product is determined uniquely by
the restriction of its norm to \(A_1\otimes A_2\).  A norm on
\(A_1\otimes A_2\) is allowed if it is a \emph{\Cstar{}norm}, that is,
multiplication and involution have norm~\(1\) and
\(\norm{x^*x}=\norm{x}^2\) for all \(x\in A_1\otimes A_2\).

There is a maximal \Cstar{}norm on \(A_1\otimes A_2\).  The resulting
\Cstar{}tensor product is called \emph{maximal \Cstar{}tensor product}
and denoted \(A_1\maxotimes A_2\).  It is characterised by the following
universal property:

\begin{proposition}
  \label{pro:tensor_max_universal}
  There is a natural bijection between non-degenerate
  \Star{}homomorphisms \(A_1\maxotimes A_2\to\Bound(\Hils)\) and pairs
  of commuting non-degenerate \Star{}homomorphisms
  \(A_1\to\Bound(\Hils)\) and \(A_2\to\Bound(\Hils)\); here we may
  replace \(\Bound(\Hils)\) by any multiplier algebra \(\Mult(D)\) of a
  \Cstar{}algebra~\(D\).
\end{proposition}

A \Star{}representation \(A\to\Bound(\Hils)\) is \emph{non-degenerate}
if \(A\cdot\Hils\) is dense in~\(\Hils\); we need this to get
representations of \(A_1\) and~\(A_2\) out of a representation of
\(A_1\maxotimes A_2\) because, for non-unital algebras, \(A_1\maxotimes
A_2\) need not contain copies of \(A_1\) and~\(A_2\).

The maximal tensor product is \emph{natural}, that is, it defines a
bifunctor
\[
\maxotimes\colon \Cstarcat\times\Cstarcat\to\Cstarcat.
\]
If \(A_1\) and~\(A_2\) are \(G\)-\Cstar{}algebras, then
\(A_1\maxotimes A_2\) inherits two group actions of~\(G\) by
naturality; these are again strongly continuous, so that
\(A_1\maxotimes A_2\) becomes a \(G\times G\)-\Cstar{}algebra.
Restricting the action to the diagonal in \(G\times G\), we turn
\(A_1\maxotimes A_2\) into a \(G\)-\Cstar{}algebra.  Thus we get a
bifunctor
\[
\maxotimes\colon \CstarcatG\times\CstarcatG\to\CstarcatG.
\]

The following lemma asserts, roughly speaking, that this tensor product
has the same formal properties as the usual tensor product for vector
spaces:

\begin{lemma}
  \label{lem:maxotimes_smc}
  There are canonical isomorphisms
  \begin{align*}
    (A\maxotimes B)\maxotimes C &\cong A\maxotimes (B\maxotimes C),\\
    A\maxotimes B &\cong B\maxotimes A,\\
    \C\maxotimes A &\cong A\cong A\maxotimes\C
  \end{align*}
  for all objects of \(\CstarcatG\) (and, in particular, of
  \(\Cstarcat\)).  These define a structure of symmetric monoidal
  category on \(\CstarcatG\) (see
  \cites{MacLane:Categories,Saavedra-Rivano:Tannakian}).
\end{lemma}

A functor between symmetric monoidal categories is called
\emph{symmetric monoidal} if it is compatible with the tensor products
in a suitable sense \cite{Saavedra-Rivano:Tannakian}.  A trivial
example is the functor \(\tau\colon \Cstarcat\to\CstarcatG\) that equips
a \Cstar{}algebra with the trivial \(G\)\nbd{}action.

It follows from the universal property that \(\maxotimes\) is compatible
with \emph{full} crossed products: if \(A\inOb\CstarcatG\),
\(B\inOb\Cstarcat\), then there is a natural isomorphism
\begin{equation}
  \label{eq:fullcross_maxotimes}
  G\cross \bigl(A\maxotimes \tau(B)\bigr) \cong (G\cross A)\maxotimes B.
\end{equation}
Like full crossed products, the maximal tensor product may be hard to
describe because it involves a maximum of all possible \Cstar{}tensor
norms.  There is another \Cstar{}tensor norm that is defined more
concretely and that combines well with \emph{reduced} crossed products.

Recall that any \Cstar{}algebra~\(A\) can be represented faithfully on a
\textsc{Hilbert} space.  That is, there is an injective
\Star{}homomorphism \(A\to\Bound(\Hils)\) for some \textsc{Hilbert}
space~\(\Hils\); here \(\Bound(\Hils)\) denotes the \Cstar{}algebra of
bounded operators on~\(\Hils\).  If~\(A\) is separable, we can find such
a representation on the separable \textsc{Hilbert} space
\(\Hils=\ell^2(\N)\).  The tensor product of two \textsc{Hilbert} spaces
\(\Hils_1\) and~\(\Hils_2\) carries a canonical inner product and can be
completed to a \textsc{Hilbert} space, which we denote by \(\Hils_1
\barotimes \Hils_2\).  If \(\Hils_1\) and~\(\Hils_2\) support faithful
representations of \Cstar{}algebras \(A_1\) and~\(A_2\), then we get an
induced \Star{}representation of \(A_1\otimes A_2\) on
\(\Hils_1\barotimes\Hils_2\).

\begin{definition}
  \label{def:spatial_tensor}
  The \emph{minimal tensor product} \(A_1 \minotimes A_2\) is the
  completion of \(A_1\otimes A_2\) with respect to the operator norm
  from \(\Bound(\Hils_1\barotimes\Hils_2)\).
\end{definition}

It can be check that this is well-defined, that is, the \Cstar{}norm on
\(A_1\otimes A_2\) does not depend on the chosen faithful
representations of \(A_1\) and~\(A_2\).  The same argument also yields
the naturality of \(A_1\minotimes A_2\).  Hence we get a bifunctor
\[ \minotimes\colon \CstarcatG\times\CstarcatG \to \CstarcatG;
\] it defines another symmetric monoidal category structure on
\(\CstarcatG\).

We may also call \(A_1\minotimes A_2\) the \emph{spatial} tensor
product.  It is \emph{minimal} in the sense that it is dominated by any
\Cstar{}tensor norm on \(A_1\otimes A_2\) that is compatible with the
given norms on \(A_1\) and~\(A_2\).  In particular, we have a canonical
surjective \Star{}homomorphism
\begin{equation}
  \label{eq:compare_max_min_otimes}
  A_1\maxotimes A_2 \to A_1\minotimes A_2.
\end{equation}

\begin{definition}
  \label{def:nuclear_Cstar}
  A \Cstar{}algebra~\(A_1\) is \emph{nuclear} if the map
  in~\eqref{eq:compare_max_min_otimes} is an isomorphism for all
  \Cstar{}algebras~\(A_2\).
\end{definition}

The name comes from an analogy between nuclear \Cstar{}algebras and
nuclear locally convex topological vector spaces
(see~\cite{Grothendieck:Produits}).  But this is merely an analogy: the
only \Cstar{}algebras that are nuclear as locally convex topological
vector spaces are the finite-dimensional ones.

Many important \Cstar{}algebras are nuclear.  This includes the
following examples:
\begin{itemize}
\item commutative \Cstar{}algebras;

\item \Cstar{}algebras of type~I and, in particular, continuous trace
\Cstar{}algebras;

\item group \Cstar{}algebras of \emph{amenable} groups (or groupoids);

\item matrix algebras and algebras of compact operators on
\textsc{Hilbert} spaces.

\end{itemize}

If~\(A\) is nuclear, then there is only one reasonable \Cstar{}algebra
completion of \(A\otimes B\).  Therefore, if we can write down any, it
must be equal to both \(A\minotimes B\) and \(A\maxotimes B\).

\begin{example}
  \label{exa:tensor_cont_X}
  For a compact space~\(X\) and a \Cstar{}algebra~\(A\), we let
  \(\CONT(X,A)\) be the \Cstar{}algebra of all continuous functions
  \(X\to A\).  If~\(X\) is a \emph{pointed} compact space, we let
  \(\CONT_0(X,A)\) be the \Cstar{}algebra of all continuous functions
  \(X\to A\) that vanish at the base point of~\(X\); this contains
  \(\CONT(X,A)\) as a special case because \(\CONT(X,A) \cong
  \CONT_0(X_+,A)\), where \(X_+=X\sqcup\{\pt\}\) with base
  point~\(\pt\).  We have
  \[
  \CONT_0(X,A)
  \cong \CONT_0(X) \minotimes A
  \cong \CONT_0(X) \maxotimes A.
  \]
\end{example}

\begin{example}
  \label{exa:tensor_compact}
  There is a unique \Cstar{}norm on \(\Mat_n\otimes A=\Mat_n(A)\) for
  all \(n\in\N\).

  For a \textsc{Hilbert} space~\(\Hils\), let \(\Comp(\Hils)\) be the
  \Cstar{}algebra of compact operators on~\(\Hils\).  Then
  \(\Comp(\Hils)\otimes A\) contains copies of \(\Mat_n(A)\),
  \(n\in\N\), for all finite-dimensional subspaces of~\(\Hils\).  These
  carry a unique \Cstar{}norm.  The \Cstar{}norms on these subspaces are
  compatible and extend to the unique \Cstar{}norm on
  \(\Comp(\Hils)\otimes A\).
\end{example}

The class of nuclear \Cstar{}algebras is closed under ideals, quotients
(by ideals), extensions, inductive limits, and crossed products by
actions of amenable locally compact groups.  In particular, this covers
crossed products by automorphisms (see Example~\ref{exa:Z_action}).

\Cstar{}subalgebras of nuclear \Cstar{}algebras need not be nuclear any
more, but they still enjoy a weaker property called \emph{exactness}:

\begin{definition}
  \label{def:Cstar_exact}
  A \Cstar{}algebra~\(A\) is called exact if the functor
  \(A\minotimes\blank\) preserves \Cstar{}algebra extensions.
\end{definition}

It is known \cites{Kirchberg-Wassermann:Exact, Ozawa:Amenable_action}
that a discrete group is exact (Definition~\ref{def:exact_group}) if and
only if its group \Cstar{}algebra is exact
(Definition~\ref{def:Cstar_exact}), if and only if the group has an
amenable action on some compact topological space.

\begin{example}
  \label{exa:free_group_exact}
  Let~\(G\) be the non-\textsc{Abel}ian free group on~\(2\) generators.
  Let~\(G\) act freely and properly on a tree as usual.  Let~\(X\) be
  the ends compactification of this tree, equipped with the induced
  action of~\(G\).  This action is known to be amenable, so that~\(G\)
  is an exact group.  Since the action is amenable, the crossed product
  algebras \(G\rcross \CONT(X)\) and \(G\cross \CONT(X)\) coincide and
  are nuclear.  The embedding \(\C\to \CONT(X)\) induces an embedding
  \(\Cred(G)\to G\cross \CONT(X)\).  But~\(G\) is not amenable.  Hence
  the \Cstar{}algebra \(\Cred(G)\) is exact but not nuclear.
\end{example}

As for crossed products, \(\minotimes\) respects injections and
surjections.  The issue with exactness in the middle is the following.
Elements of \(A\minotimes B\) are limits of tensors of the form
\(\sum_{i=1}^n a_i\otimes b_i\) with \(a_1,\dotsc,a_n\in A\),
\(b_1,\dotsc,b_n\in B\).  If an element in \(A\minotimes B\) is
annihilated by the map to \((A/I)\minotimes B\), then we can approximate
it by such finite sums for which \(\sum_{i=1}^n (a_i\bmod I)\otimes
b_i\) goes to~\(0\).  But this does not suffice to find approximations
in \(I\otimes B\).  Thus the kernel of the projection map \(A\minotimes
B\prto (A/I)\minotimes B\) may be strictly larger than \(I\minotimes
B\).

\begin{proposition}
  \label{pro:exact_on_cpsplit_tensor}
  The functor \(\maxotimes\) is exact in each variable, that is,
  \(\blank\maxotimes D\) maps extensions in \(\CstarcatG\) again to
  extensions for each \(D\inOb\CstarcatG\).

  Both \(\minotimes\) and \(\maxotimes\) map split extensions to split
  extensions and (equivariantly) cp-split extensions again to
  (equivariantly) cp-split extensions.
\end{proposition}

The proof is similar to the proofs of Propositions
\ref{pro:fullcross_exact} and~\ref{pro:exact_on_cpsplit}.

If \(A\) or~\(B\) is nuclear, we simply write \(A\otimes B\) for
\(A\maxotimes B\cong A\minotimes B\).

\subsection{Limits and colimits}
\label{sec:limits}

\begin{proposition}
  \label{pro:limits}
  The categories \(\Cstarcat\) and \(\CstarcatG\) are bicomplete, that
  is, any (small) diagram in these categories has both a limit and a
  colimit.
\end{proposition}

\begin{proof}
  To get general limits and colimits, it suffices to construct
  equalisers and coequalisers for pairs of parallel morphisms
  \(f_0,f_1\colon A\rightrightarrows B\), direct products and coproducts
  \(A_1\times A_2\) and \(A_1\sqcup A_2\) for any pair of objects and,
  more generally, for arbitrary sets of objects.

  The equaliser and coequaliser of \(f_0,f_1\colon A\rightrightarrows
  B\) are
  \[
  \ker(f_0-f_1) = \{a\in A\mid f_0(a)=f_1(a)\}\subseteq A
  \]
  and the quotient of~\(A_1\) by the closed \Star{}ideal generated by
  the range of \(f_0-f_1\), respectively.  Here we use that quotients of
  \Cstar{}algebras by closed \Star{}ideals are again \Cstar{}algebras.
  Notice that \(\ker(f_0-f_1)\) is indeed a \Cstar{}subalgebra of~\(A\).

  The direct product \(A_1\times A_2\) is the usual direct product,
  equipped with the canonical \Cstar{}algebra structure.  We can
  generalise the construction of the direct product to infinite direct
  products: let \(\prod_{i\in I} A_i\) be the set of all
  \emph{norm-bounded} sequences \((a_i)_{i\in I}\) with \(a_i\in A_i\)
  for all \(i\in I\); this is a \Cstar{}algebra with respect to the
  obvious \Star{}algebra structure and the norm \(\norm{(a_i)} \defeq
  \sup_{i\in I} \norm{a_i}\).  It has the right universal property
  because any \Star{}homomorphism is norm-contracting.  (A similar
  construction with \textsc{Banach} algebras would fail at this point.)

  The coproduct \(A_1\sqcup A_2\) is also called \emph{free product} and
  denoted \(A_1*A_2\); its construction is more involved.  The free
  \(\C\)\nbd{}algebra generated by \(A_1\) and~\(A_2\) carries a
  canonical involution, so that it makes sense to study \Cstar{}norms on
  it.  It turns out that there is a maximal such \Cstar{}norm.  The
  resulting \Cstar{}completion is the free product \Cstar{}algebra.  In
  the equivariant case, \(A_1\sqcup A_2\) inherits an action of~\(G\),
  which is strongly continuous.  The resulting object of \(\CstarcatG\)
  has the correct universal property for a coproduct.

  An \emph{inductive system} of \Cstar{}algebras
  \((A_i,\alpha_i^j)_{i\in I}\) is called \emph{reduced} if all the maps
  \(\alpha_i^j\colon A_i\to A_j\) are injective; then they are
  automatically isometric embeddings.  We may as well assume that these
  maps are identical inclusions of \Cstar{}subalgebras.  Then we can
  form a \Star{}algebra \(\bigcup A_i\), and the given \Cstar{}norms
  piece together to a \Cstar{}norm on \(\bigcup A_i\).  The resulting
  completion is \(\varinjlim (A_i,\alpha_i^j)\).  In particular, we can
  construct an infinite coproduct as the inductive limit of its finite
  sub-coproducts.  Thus we get infinite coproducts.
\end{proof}

The category of commutative \Cstar{}algebras is equivalent to the
opposite of the category of pointed compact spaces by the
\textsc{Gelfand}--\textsc{Naimark} Theorem.  It is frequently convenient
to replace a \emph{pointed compact} space~\(X\) with base point~\(\pt\)
by the \emph{locally compact} space \(X\setminus\{\pt\}\).  A continuous
map \(X\to Y\) extends to a pointed continuous map \(X_+\to Y_+\) if and
only if it is \emph{proper}.  But there are more pointed continuous maps
\(f\colon X_+\to Y_+\) than proper continuous maps \(X\to Y\) because
points in~\(X\) may be mapped to the point at infinity \(\infty\in
Y_+\).  For instance, the zero homomorphism \(\CONT_0(Y)\to\CONT_0(X)\)
corresponds to the constant map \(x\mapsto\infty\).

\begin{example}
  \label{exa:open_subset_embed}
  If \(U\subseteq X\) is an open subset of a locally compact space, then
  \(\CONT_0(U)\) is an ideal in \(\CONT_0(X)\).  No map \(X\to U\)
  corresponds to the embedding \(\CONT_0(U)\to\CONT_0(X)\).
\end{example}

\begin{example}
  \label{exa:infinite_product_bad}
  Products of commutative \Cstar{}algebras are again commutative and
  correspond by the \textsc{Gelfand}--\textsc{Naimark} Theorem to
  coproducts in the category of pointed compact spaces.  The coproduct
  of a set of pointed compact spaces is the
  \emph{\textsc{Stone}--\textsc{\v{C}ech} compactification} of their
  wedge sum.  Thus infinite products in \(\Cstarcat\) and~\(\CstarcatG\)
  do not behave well for the purposes of homotopy theory.

  The coproduct of two non-zero \Cstar{}algebras is never commutative
  and hence has no analogue for (pointed) compact spaces.  The smash
  product for pointed compact spaces corresponds to the tensor product
  of \Cstar{}algebras because
  \[
  \CONT_0(X\wedge Y) \cong \CONT_0(X)\minotimes \CONT_0(Y).
  \]
\end{example}

\subsection{Enrichment over pointed topological spaces}
\label{sec:enrichment}

Let \(A\) and~\(B\) be \Cstar{}algebras.  It is well-known that a
\Star{}homomorphism \(f\colon A\to B\) is automatically norm-contracting
and induces an isometric embedding \(A/\ker f \to B\) with respect to
the quotient norm on \(A/\ker f\).  The reason for this is that the norm
for self-adjoint elements in a \Cstar{}algebra agrees with the spectral
radius and hence is determined by the algebraic structure; by the
\Cstar{}condition \(\norm{a}^2=\norm{a^*a}\), this extends to all
elements of a \Cstar{}algebra.

It follows that \(\Hom(A,B)\) is an \emph{equicontinuous} set of
linear maps \(A\to B\).  We always equip \(\Hom(A,B)\) with the
topology of pointwise norm-convergence.  Its subbasic open subsets are
of the form
\[
\{f\colon A\to B\mid \norm{(f-f_0)(a)}<1 \ \forall a\in S\}
\]
for \(f_0\in\Hom(A,B)\) and a finite subset \(S\subseteq A\).  Since
\(\Hom(A,B)\) is equicontinuous, this topology agrees with the
topology of uniform convergence on compact subsets, which is generated
by the corresponding subsets for compact~\(S\).  This is nothing but
the compact-open topology on mapping spaces.  But it differs from the
topology defined by the operator norm.  We shall never use the latter.

\begin{lemma}
  \label{lem:Hom_metrisable}
  If~\(A\) is separable, then \(\Hom(A,B)\) is metrisable for any~\(B\).
\end{lemma}

\begin{proof}
  There exists a sequence \((a_n)_{n\in\N}\) in~\(A\) with \(\lim
  a_n=0\) whose closed linear span is all of~\(A\).  The metric
  \[
  d(f_1,f_2) = \sup \{\norm{f_1(a_n) - f_2(a_n)} \mid n\in\N\}
  \]
  defines the topology of uniform convergence on compact subsets on
  \(\Hom(A,B)\) because the latter is equicontinuous.
\end{proof}

There is a distinguished element in \(\Hom(A,B)\) as well, namely, the
\emph{zero homomorphism} \(A\to 0\to B\).  Thus \(\Hom(A,B)\) becomes a
pointed topological space.

\begin{proposition}
  \label{pro:Cstar_enriched_top}
  The above construction provides an enrichment of \(\Cstarcat\) over
  the category of pointed Hausdorff topological spaces.
\end{proposition}

\begin{proof} It is clear that \(0\circ f = 0\) and \(f\circ 0 = 0\) for
all morphisms~\(f\).  Furthermore, we must check that composition of
morphisms is \emph{jointly} continuous.  This follows from the
equicontinuity of \(\Hom(A,B)\).
\end{proof}

This enrichment allows us to carry over some important definitions from
categories of spaces to \(\Cstarcat\).  For instance, a \emph{homotopy}
between two \Star{}homomorphisms \(f_0,f_1\colon A\to B\) is a
continuous path between \(f_0\) and~\(f_1\) in the topological space
\(\Hom(A,B)\).  In the following proposition, \(\Map_+(X,Y)\) denotes
the space of morphisms in the category of pointed topological spaces,
equipped with the compact-open topology.

\begin{proposition}[compare Proposition 3.4 in
  \cite{Joachim-Johnson:Model}]
  \label{pro:Map_Hom}
  Let \(A\) and~\(B\) be \Cstar{}algebras and let~\(X\) be a pointed
  compact space.  Then
  \[
  \Map_+\bigl(X,\Hom(A,B)\bigr) \cong \Hom\bigl(A,\CONT_0(X,B)\bigr)
  \]
  as pointed topological spaces.
\end{proposition}

\begin{proof}
  If we view \(A\) and~\(B\) as pointed topological spaces with the norm
  topology and base point~\(0\), then \(\Hom(A,B)\subseteq \Map_+(A,B)\)
  is a topological subspace with the same base point because
  \(\Hom(A,B)\) also carries the compact-open topology.  Since~\(X\) is
  compact, standard point set topology yields homeomorphisms
  \begin{multline*}
    \Map_+\bigl(X,\Map_+(A,B)\bigr) \cong \Map_+(X\wedge A,B)
    \\\cong \Map_+\bigl(A,\Map_+(X,B)\bigr)
    = \Map_+\bigl(A,\CONT_0(X,B)\bigr).
  \end{multline*}
  These restrict to the desired homeomorphism.
\end{proof}

In particular, a homotopy between two \Star{}homomorphisms
\(f_0,f_1\colon A\rightrightarrows B\) is equivalent to a
\Star{}homomorphism \(f\colon A\to \CONT([0,1],B)\) with \(\ev_t\circ
f=f_t\) for \(t=0,1\), where \(\ev_t\) denotes the \Star{}homomorphism
\[
\ev_t\colon \CONT([0,1],B)\to B, \qquad f\mapsto f(t).
\]

We also have
\[
\CONT_0\bigl(X,\CONT_0(Y,A)\bigr) \cong \CONT_0(X\wedge Y,A)
\]
for all pointed compact spaces \(X\), \(Y\) and all
\(G\)-\Cstar{}algebras~\(A\).  Thus a homotopy between two homotopies
can be encoded by a \Star{}homomorphism
\[
A\to C\bigl([0,1],\CONT([0,1],B)\bigr) \cong \CONT([0,1]^2,B).
\]

These constructions work only for pointed \emph{compact} spaces.  If we
enlarge the category of \Cstar{}algebras to a suitable category of
projective limits of \Cstar{}algebras as
in~\cite{Joachim-Johnson:Model}, then we can define \(\CONT_0(X,A)\) for
any pointed \emph{compactly generated} space~\(X\).  But we lose some of
the nice analytic properties of \Cstar{}algebras.  Therefore, I prefer
to stick to the category of \Cstar{}algebras itself.

\subsection{Cylinders, cones, and suspensions}
\label{sec:cylinders}

The following definitions go back to \cite{Schochet:Top3}, where some
more results can be found.  The description of homotopies above leads us
to define the \emph{cylinder} over a \Cstar{}algebra~\(A\) by
\[
\Cyl(A) \defeq \CONT([0,1],A).
\]
This is compatible with the cylinder construction for spaces because
\[
\Cyl\bigl(\CONT_0(X)\bigr) \cong \CONT\bigl([0,1],\CONT_0(X)\bigr) \cong
\CONT_0([0,1]_+\wedge X)
\]
for any pointed compact space~\(X\); if we use locally compact spaces,
we get \([0,1]\times X\) instead of \([0,1]_+\wedge X\).

The universal property of \(\Cyl(A)\) is dual to the usual one for spaces
because the identification between pointed compact spaces and commutative
\Cstar{}algebras is contravariant.

Similarly, we may define the \emph{cone} \(\Cone(A)\) and the
\emph{suspension} \(\Sus(A)\) by
\[
\Cone(A) \defeq  \CONT_0\bigl([0,1]\setminus\{0\},A),\qquad
\Sus(A) \defeq  \CONT_0\bigl([0,1]\setminus\{0,1\},A) \cong
\CONT_0(\Sphere^1,A),
\]
where~\(\Sphere^1\) denotes the pointed \(1\)\nbd{}sphere, that is,
circle.  These constructions are compatible with the corresponding ones
for spaces as well, that is,
\[
\Cone\bigl(\CONT_0(X)\bigr) \cong \CONT_0([0,1]\wedge X),\quad
\Sus\bigl(\CONT_0(X)\bigr) \cong \CONT_0( \Sphere^1\wedge X).
\]
Here \([0,1]\) has the base point~\(0\).

\begin{definition}
  \label{def:mapping_cone}
  Let \(f\colon A\to B\) be a morphism in \(\Cstarcat\) or
  \(\CstarcatG\).  The \emph{mapping cylinder} \(\Cyl(f)\) and the
  \emph{mapping cone} \(\Cone(f)\) of~\(f\) are the limits of the
  diagrams
  \[
  A\xrightarrow{f} B \xleftarrow{\ev_1} \Cyl(B),
  \qquad
  A\xrightarrow{f} B \xleftarrow{\ev_1} \Cone(B).
  \]
\end{definition}

More concretely,
\begin{align*}
  \Cone(f) &= \bigl\{(a,b)\in A\times \CONT_0\bigl((0,1],B\bigr)
  \bigm| f(a)= b(1)\bigr\},\\
  \Cyl(f) &= \bigl\{(a,b)\in A\times \CONT_0([0,1],B\bigr)
  \bigm| f(a)= b(1)\bigr\}.
\end{align*}

If \(f\colon X\to Y\) is a morphism of pointed compact spaces, then the
mapping cone and mapping cylinder of the induced \Star{}homomorphism
\(\CONT_0(f) \colon \CONT_0(Y)\to\CONT_0(X)\) agree with
\(\CONT_0\bigl(\Cyl(f)\bigr)\) and \(\CONT_0\bigl(\Cone(f)\bigr)\),
respectively.

The cylinder, cone, and suspension functors are exact for various kinds
of extensions: they map extensions, split extensions, and cp-split
extensions again to extensions, split extensions, and cp-split
extensions, respectively.  Similar remarks apply to mapping cylinders
and mapping cones: for any morphism of extensions
\begin{equation}
  \label{eq:morphism_of_extensions}
  \begin{gathered}
    \xymatrix{
      I\ \ar@{>->}[r] \ar[d]^{\alpha}&
      E \ar@{->>}[r] \ar[d]^{\beta}&
      Q \ar[d]^{\gamma}\\
      I'\ \ar@{>->}[r] &
      E' \ar@{->>}[r] &
      Q', \\
    }
  \end{gathered}
\end{equation}
we get extensions
\begin{equation}
  \label{eq:cyl_cone_exact}
  \Cyl(\alpha)\into \Cyl(\beta)\prto \Cyl(\gamma),
  \qquad
  \Cone(\alpha)\into \Cone(\beta)\prto \Cone(\gamma);
\end{equation}
if the extensions in~\eqref{eq:morphism_of_extensions} are split or
cp-split, so are the resulting extensions in~\eqref{eq:cyl_cone_exact}.

The familiar maps relating mapping cones and cylinders to cones and
suspensions continue to exist in our case.  For any morphism \(f\colon
A\to B\) in \(\CstarcatG\), we get a morphism of extensions
\[\xymatrix{
  \Sus(B)\ \ar@{>->}[r] \ar[d] & \Cone(f) \ar@{->>}[r] \ar[d] &
  A \ar@{=}[d]\\
  \Cone(B)\ \ar@{>->}[r] & \Cyl(f) \ar@{->>}[r] & A\\
}
\]
The bottom extension splits and the maps \(A\leftrightarrow \Cyl(f)\)
are inverse to each other up to homotopy.  By naturality, the composite
map \(\Cone(f) \to A\to B\) factors through \(\Cone(\ID_B) \cong
\Cone(B)\) and hence is homotopic to the zero map.

\section{Universal functors with certain properties}
\label{sec:properties_functors}

When we study topological invariants for \Cstar{}algebras, we usually
require homotopy invariance and some exactness and stability conditions.
Here we investigate these conditions and their interplay and describe
some universal functors.

We discuss homotopy invariant functors on \(\CstarcatG\) and the
homotopy category \(\Ho(\CstarcatG)\) in \textsection\ref{sec:homotopy}.
This is parallel to classical topology.  We turn to
\textsc{Morita}--\textsc{Rieffel} invariance and \Cstar{}stability in
\textsection\ref{sec:MR_equivalence}--\ref{sec:Cstar_stable}.  We
describe the resulting localisation using correspondences.  By the way,
in a \Cstar{}algebra context, \emph{stability} usually refers to
algebras of compact operators instead of suspensions.
\textsection\ref{sec:exactness} deals with various exactness conditions:
split-exactness, half-exactness, and additivity.

Whereas each of the above properties in itself seems rather weak, their
combination may have striking consequences.  For instance, a functor
that is both \Cstar{}stable and split-exact is automatically homotopy
invariant and satisfies \textsc{Bott} periodicity.

Throughout this section, we consider functors \(\Cstarnice\to\Cat\)
where~\(\Cstarnice\) is a full subcategory of \(\Cstarcat\) or
\(\CstarcatG\) for some locally compact group~\(G\).  The target
category~\(\Cat\) may be arbitrary in
\textsection\ref{sec:homotopy}--\textsection\ref{sec:Cstar_stable}; to
discuss exactness properties, we require~\(\Cat\) to be an exact
category or at least additive.  Typical choices for~\(\Cstarnice\) are
the categories of separable or separable nuclear \Cstar{}algebras, or
the subcategory of all separable nuclear \(G\)\nbd{}\Cstar{}algebras
with an amenable (or a proper) action of~\(G\).

\begin{definition}
  \label{def:univ_functor}
  Let~\(P\) be a property for functors defined on~\(\Cstarnice\).  A
  \emph{universal functor} with~\(P\) is a functor \(u\colon
  \Cstarnice\to \Univ_P(\Cstarnice)\) such that
  \begin{itemize}
  \item \(\bar{F}\circ u\) has~\(P\) for each functor \(\bar{F}\colon
    \Univ_P(\Cstarnice)\to\Cat\);

  \item any functor \(F\colon \Cstarnice\to\Cat\) with~\(P\) factors
    uniquely as \(F=\bar{F}\circ u\) for some functor \(\bar{F}\colon
    \Univ_P(\Cstarnice)\to\Cat\).
  \end{itemize}
\end{definition}

Of course, a universal functor with~\(P\) need not exist.  If it does,
then it restricts to a bijection between objects of~\(\Cstarnice\) and
\(\Univ_P(\Cstarnice)\).  Hence we can completely describe it by the
sets of morphisms \(\Univ_P(A,B)\) from~\(A\) to~\(B\) in
\(\Univ_P(\Cstarnice)\) and the maps \(\Cstarnice(A,B)\to\Univ_P(A,B)\)
for \(A,B\inOb\Cstarnice\); the universal property means that for any
functor \(F\colon \Cstarnice\to\Cat\) with~\(P\) there is a unique
functorial way to extend the maps \(\Hom_G(A,B)\to
\Cat\bigl(F(A),F(B)\bigr)\) to \(\Univ_P(A,B)\).  There is no \textit{a
  priori} reason why the morphism spaces \(\Univ_P(A,B)\) for
\(A,B\inOb\Cstarnice\) should be independent of~\(\Cstarnice\); but this
happens in the cases we consider, under some assumption
on~\(\Cstarnice\).

\subsection{Homotopy invariance}
\label{sec:homotopy}

The following discussion applies to any full subcategory
\(\Cstarnice\subseteq\CstarcatG\) that is closed under the cylinder,
cone, and suspension functors.

\begin{definition}
  \label{def:homotopy_category}
  Let \(f_0,f_1\colon A\rightrightarrows B\) be two parallel morphisms
  in~\(\Cstarnice\).  We write \(f_0\sim f_1\) and call \(f_0\)
  and~\(f_1\) \emph{homotopic} if there is a \emph{homotopy} between
  \(f_0\) and~\(f_1\), that is, a morphism \(f\colon A\to \Cyl(B) =
  \CONT([0,1],B)\) with \(\ev_t\circ f=f_t\) for \(t=0,1\).

  It is easy to check that homotopy is an equivalence relation on
  \(\Hom_G(A,B)\).  We let \([A,B]\) be the set of equivalence classes.
  The composition of morphisms in~\(\Cstarnice\) descends to maps
  \[
  [B,C]\times [A,B]\to [A,C],
  \qquad \bigl([f],[g]\bigr)\mapsto [f\circ g],
  \]
  that is, \(f_1\sim f_2\) and \(g_1\sim g_2\) implies \(f_1\circ
  f_2\sim g_1\circ g_2\).  Thus the sets \([A,B]\) form the morphism
  sets of a category, called \emph{homotopy category of~\(\Cstarnice\)}
  and denoted \(\Ho(\Cstarnice)\).  The identity maps on objects and the
  canonical maps on morphisms define a \emph{canonical functor}
  \(\Cstarnice\to\Ho(\Cstarnice)\).  A morphism in~\(\Cstarnice\) is
  called a \emph{homotopy equivalence} if it becomes invertible in
  \(\Ho(\Cstarnice)\).
\end{definition}

\begin{lemma}
  \label{lem:homotopy_category}
  The following are equivalent for a functor \(F\colon
  \Cstarnice\to\Cat\):
  \begin{enumerate}[(a)]
  \item \(F(\ev_0)=F(\ev_1)\) as maps \(F\bigl(\CONT([0,1],A)\bigr)\to
    F(A)\) for all \(A\inOb\Cstarnice\);

  \item \(F(\ev_t)\) induces isomorphisms
    \(F\bigl(\CONT([0,1],A)\bigr)\to F(A)\) for all
    \(A\inOb\Cstarnice\), \(t\in[0,1]\);

  \item the embedding as constant functions \(\const\colon A\to
    \CONT([0,1],A)\) induces an isomorphism \(F(A) \to
    F\bigl(\CONT([0,1],A)\bigr)\) for all \(A\inOb\Cstarnice\);

  \item \(F\) maps homotopy equivalences to isomorphisms;

  \item if two parallel morphisms \(f_0,f_1\colon A\rightrightarrows B\)
    are homotopic, then \(F(f_0)=F(f_1)\);

  \item \(F\) factors through the canonical functor
    \(\Cstarnice\to\Ho(\Cstarnice)\).

  \end{enumerate}
  Furthermore, the factorisation \(\Ho(\Cstarnice)\to \Cat\) in~(f) is
  necessarily unique.
\end{lemma}

\begin{proof} We only mention two facts that are needed for the proof.
  First, we have \(\ev_t\circ\const=\ID_A\) and \(\const\circ\ev_t\sim
  \ID_{\CONT([0,1],A)}\) for all \(A\inOb\Cstarnice\) and all
  \(t\in[0,1]\).  Secondly, an isomorphism has a unique left and a
  unique right inverse, and these are again isomorphisms.
\end{proof}

The equivalence of (d) and~(f) in Lemma \ref{lem:homotopy_category} says
that \(\Ho(\Cstarnice)\) is the localisation of~\(\Cstarnice\) at the
family of homotopy equivalences.

Since \(\CONT\bigl([0,1],\CONT_0(X)\bigr) \cong \CONT_0([0,1]\times X)\)
for any locally compact space~\(X\), our notion of homotopy restricts to
the usual one for pointed compact spaces.  Hence the opposite of the
homotopy category of pointed compact spaces is equivalent to a full
subcategory of \(\Ho(\Cstarcat)\).

The sets \([A,B]\) inherit a base point \([0]\) and a quotient topology
from \(\Hom(A,B)\); thus \(\Ho(\Cstarnice)\) is enriched over pointed
topological spaces as well.  This topology on \([A,B]\) is not so
useful, however, because it need not be Hausdorff.

A similar topology exists on \textsc{Kasparov} groups and can be defined
in various ways, which turn out to be equivalent
\cite{Dadarlat:Topology}.

Let \(F\colon \CstarcatG\to\CstarcatG[H]\) be a functor with natural
isomorphisms
\[
F\bigl(\CONT([0,1],A)\bigr) \cong \CONT([0,1],F(A)\bigr)
\]
that are compatible with evaluation maps for all~\(A\).  The universal
property implies that~\(F\) descends to a functor
\(\Ho(\CstarcatG)\to\Ho(\CstarcatG[H])\).  In particular, this applies
to the suspension, cone, and cylinder functors and, more generally, to
the functors \(A\maxotimes\blank\) and \(A\minotimes\blank\) on
\(\CstarcatG\) for any \(A\inOb\CstarcatG\) because both tensor product
functors are associative and commutative and
\[
\CONT([0,1],A) \cong \CONT([0,1])\maxotimes A
\cong \CONT([0,1])\minotimes A.
\]
The same works for the reduced and full crossed product functors
\(\CstarcatG\to\Cstarcat\).

We may stabilise the homotopy category with respect to the
\emph{suspension} functor and consider a \emph{suspension-stable
  homotopy category} with morphism spaces
\[
\{A,B\} \defeq \varinjlim_{k\to\infty} [\Sus^k{A},\Sus^k{B}]
\]
for all \(A,B\inOb\CstarcatG\).  We may also enlarge the set of objects
by adding formal desuspensions and generalising the notion of
\emph{spectrum}.  This is less interesting for \Cstar{}algebras than for
spaces because most functors of interest satisfy \textsc{Bott}
Periodicity, so that suspension and desuspension become equivalent.

\subsection{\textsc{Morita}--\textsc{Rieffel} equivalence and stable
  isomorphism}
\label{sec:MR_equivalence}

One of the basic ideas of non-commutative geometry is that \(G\cross
\CONT_0(X)\) (or \(G\rcross \CONT_0(X)\)) should be a substitute for the
quotient space \(G\backslash X\), which may have bad singularities.  In
the special case of a free and proper \(G\)\nbd{}space~\(X\), we expect
that \(G\cross \CONT_0(X)\) and \(\CONT_0(G\backslash X)\) are
``equivalent'' in a suitable sense.  Already the simplest possible case
\(X=G\) shows that we cannot expect an isomorphism here because
\[
G\cross \CONT_0(G) \cong G\rcross \CONT_0(G) \cong \Comp(L^2G).
\]
The right notion of equivalence is a \Cstar{}version of \textsc{Morita}
equivalence due to \textsc{Marc A.\ Rieffel} (\cites{Rieffel:Induced,
  Rieffel:Morita, Rieffel:Strong_Morita}); therefore, we call it
\textsc{Morita}--\textsc{Rieffel} equivalence.

The definition of \textsc{Morita}--\textsc{Rieffel} equivalence involves
\textsc{Hilbert} modules over \(C^*\)\nbd{}algebras and the
\Cstar{}algebras of compact operators on them; these notions are crucial
for \textsc{Kasparov} theory as well.  We refer to
\cite{Lance:Hilbert_modules} for the definition and a discussion of
their basic properties.

\begin{definition}
  \label{def:MR_equivalence}
  Two \(G\)\nbd{}\Cstar{}algebras \(A\) and~\(B\) are called
  \emph{\textsc{Morita}--\textsc{Rieffel} equivalent} if there are a
  full \(G\)\nbd{}equivariant \textsc{Hilbert}
  \(B\)\nbd{}module~\(\Hilm\) and a \(G\)\nbd{}equivariant
  \Star{}isomorphism \(\Comp(\Hilm) \cong A\).
\end{definition}

It is possible (and desirable) to express this definition more
symmetrically: \(\Hilm\) is an \(A,B\)\brd{}bimodule with two inner
products taking values in \(A\) and~\(B\), satisfying various conditions
\cite{Rieffel:Induced}.  \textsc{Morita}--\textsc{Rieffel} equivalent
\(G\)\nbd{}\Cstar{}algebras have equivalent categories of
\(G\)\nbd{}equivariant \textsc{Hilbert} modules via
\(\Hilm\otimes_B\blank\).  The converse is unclear.

\begin{example}
  \label{exa:double_coset_MR_equivalence}
  The following is a more intricate example of a
  \textsc{Morita}--\textsc{Rieffel} equivalence.  Let \(\Gamma\)
  and~\(P\) be two subgroups of a locally compact group~\(G\).
  Then~\(\Gamma\) acts on \(G/P\) by left translation and~\(P\) acts on
  \(\Gamma\backslash G\) by right translation.  The corresponding orbit
  space is the double coset space \(\Gamma\backslash G/P\).  Both
  \(\Gamma \cross \CONT_0(G/P)\) and \(P\cross \CONT_0(\Gamma\backslash
  G)\) are non-commutative models for this double coset space.  They are
  indeed \textsc{Morita}--\textsc{Rieffel} equivalent; the bimodule that
  implements the equivalence is a suitable completion of \(\CONT_c(G)\),
  the space of continuous functions with compact support on~\(G\).
\end{example}

These examples suggest that \textsc{Morita}--\textsc{Rieffel} equivalent
\Cstar{}algebras describe the \emph{same} non-commutative space.
Therefore, we expect that reasonable functors on \(\Cstarcat\) should
not distinguish between \textsc{Morita}--\textsc{Rieffel} equivalent
\Cstar{}algebras.  (We will slightly weaken this statement below.)

\begin{definition}
  \label{def:stable_equivalence}
  Two \(G\)\nbd{}\Cstar{}algebras \(A\) and~\(B\) are called
  \emph{stably isomorphic} if there is a \(G\)\nbd{}equivariant
  \Star{}isomorphism \(A\otimes \Comp(\Hils_G) \cong B\otimes
  \Comp(\Hils_G)\), where \(\Hils_G\defeq L^2(G\times\N)\) is the direct
  sum of countably many copies of the regular representation of~\(G\);
  we let~\(G\) act on \(\Comp(\Hils_G)\) by conjugation, of course.
\end{definition}

The following technical condition is often needed in connection with
\textsc{Morita}--\textsc{Rieffel} equivalence.

\begin{definition}
  \label{def:sigma_unital}
  A \Cstar{}algebra is called \emph{\(\sigma\)\nbd{}unital} if it has a
  countable approximate identity or, equivalently, contains a strictly
  positive element.
\end{definition}

\begin{example}
  \label{exa:sigma_unital}
  All separable \Cstar{}algebras and all unital \Cstar{}algebras are
  \(\sigma\)\nbd{}unital; the algebra \(\Comp(\Hils)\) is
  \(\sigma\)\nbd{}unital if and only if~\(\Hils\) is separable.
\end{example}

\begin{theorem}[\cite{Brown-Green-Rieffel}]
  \label{the:MR_stable_equivalence}
  \(\sigma\)\nbd{}Unital \(G\)\nbd{}\Cstar{}algebras are
  \(G\)\nbd{}equivariantly \textsc{Morita}--\textsc{Rieffel} equivalent
  if and only if they are stably isomorphic.
\end{theorem}

In the non-equivariant case, this theorem is due to
\textsc{Brown}--\textsc{Green}--\textsc{Rieffel}
\cite{Brown-Green-Rieffel}.  A simpler proof that carries over to the
equivariant case appeared in \cite{Mingo-Phillips}.

\subsection{\texorpdfstring{$C^*$-}{C*-}stable functors}
\label{sec:Cstar_stable}

The definition of \Cstar{}stability is more intuitive in the
non-equivariant case:

\begin{definition}
  \label{def:Cstar_stable}
  Fix a rank-one projection \(p\in\Comp(\ell^2\N)\).  The resulting
  embedding \(A\to A\otimes\Comp(\ell^2\N)\), \(a\mapsto a\otimes p\),
  is called a \emph{corner embedding} of~\(A\).

  A functor \(F\colon \Cstarcat\to\Cat\) is called \emph{\Cstar{}stable}
  if any corner embedding induces an isomorphism \(F(A)\cong
  F\bigl(A\otimes \Comp(\ell^2\N)\bigr)\).
\end{definition}

The correct equivariant generalisation is the following:

\begin{definition}[\cite{Meyer:KKG}]
  \label{def:Cstar_stable_G}
  A functor \(F\colon \CstarcatG\to\Cat\) is called
  \emph{\Cstar{}stable} if the canonical embeddings
  \(\Hils_1\to\Hils_1\oplus\Hils_2\leftarrow \Hils_2\) induce
  isomorphisms
  \[
  F\bigl(A\otimes \Comp(\Hils_1)\bigr) \xrightarrow{\cong}
  F\bigl(A\otimes \Comp(\Hils_1\oplus \Hils_2)\bigr) \xleftarrow{\cong}
  F\bigl(A\otimes \Comp(\Hils_2)\bigr)
  \]
  for all non-zero \(G\)\nbd{}\textsc{Hilbert} spaces \(\Hils_1\)
  and~\(\Hils_2\).
\end{definition}

Of course, it suffices to require \(F\bigl(A\otimes \Comp(\Hils_1)\bigr)
\xrightarrow{\cong} F\bigl(A\otimes \Comp(\Hils_1\oplus
\Hils_2)\bigr)\).  It is not hard to check that Definitions
\ref{def:Cstar_stable} and~\ref{def:Cstar_stable_G} are equivalent for
trivial~\(G\).

\begin{remark}
  \label{rem:matrix_stable}
  We have argued in \textsection\ref{sec:MR_equivalence} why
  \Cstar{}stability is an essential property for any decent homology
  theory for \Cstar{}algebras.  Nevertheless, it is tempting to assume
  less because \Cstar{}stability together with split-exactness has very
  strong implications.

  One reasonable way to weaken \Cstar{}stability is to replace
  \(\Comp(\ell^2\N)\) by~\(\Mat_n\) for \(n\in\N\) in
  Definition~\ref{def:Cstar_stable} (see~\cite{Thom:Connective}).  If
  two \emph{unital} \Cstar{}algebras are
  \textsc{Morita}--\textsc{Rieffel} equivalent, then they are also
  \textsc{Morita} equivalent as rings, that is, the equivalence is
  implemented by a \emph{finitely generated} projective module.  This
  implies that a matrix-stable functor is invariant under
  \textsc{Morita}--\textsc{Rieffel} equivalence for \emph{unital}
  \Cstar{}algebras.

  Matrix-stability also makes good sense in \(\CstarcatG\) for a
  \emph{compact} group~\(G\): simply require \(\Hils_1\) and~\(\Hils_2\)
  in Definition~\ref{def:Cstar_stable_G} to be finite-dimensional.  But
  we seem to run into problems for non-compact groups because they may
  have few finite-dimensional representations and we lack a
  finite-dimensional version of the equivariant stabilisation theorem.
\end{remark}

Our next goal is to describe the \emph{universal} \Cstar{}stable
functor.  We abbreviate \(A_\Comp\defeq \Comp(L^2G)\otimes A\).

\begin{definition}
  \label{def:correspondence}
  A \emph{correspondence} from~\(A\) to~\(B\) (or \(A\dasharrow B\)) is
  a \(G\)\nbd{}equivariant \textsc{Hilbert}
  \(B_\Comp\)\nbd{}module~\(\Hilm\) together with a
  \(G\)\nbd{}equivariant essential (or non-degenerate)
  \Star{}homomorphism \(f\colon A_\Comp\to\Comp(\Hilm)\).

  Given correspondences~\(\Hilm\) from~\(A\) to~\(B\) and~\(\Hilm[F]\)
  from~\(B\) to~\(C\), their \emph{composition} is the correspondence
  from~\(A\) to~\(C\) with underlying \textsc{Hilbert} module
  \(\Hilm\barotimes_{B_\Comp} \Hilm[F]\) and map \(A_\Comp\to
  \Comp(\Hilm) \to \Comp(\Hilm\barotimes_{B_\Comp} \Hilm[F])\), where
  the last map sends \(T\mapsto T\otimes 1\); this yields compact
  operators because~\(B_\Comp\) maps to \(\Comp(\Hilm[F])\).
  See~\cite{Lance:Hilbert_modules} for the definition of the relevant
  completed tensor product of \textsc{Hilbert} modules.
\end{definition}

Up to isomorphism, the composition of correspondences is associative and
the identity maps \(A\to A=\Comp(A)\) act as unit elements.  Hence we
get a category \(\Corr{G}\) whose morphisms are the \emph{isomorphism
  classes} of correspondences.  It may have advantages to treat
\(\Corr{G}\) as a \(2\)\nbd{}category.

Any \Star{}homomorphism \(\varphi\colon A\to B\) yields a correspondence
\(f\colon A\to \Comp(\Hilm)\) from \(A\) to~\(B\), so that we get a
canonical functor \(\natural\colon \CstarcatG\to\Corr{G}\).  We
let~\(\Hilm\) be the right ideal \(\varphi(A_\Comp)\cdot B_\Comp\)
in~\(B_\Comp\), viewed as a \textsc{Hilbert} \(B\)\nbd{}module.  Then
\(f(a)\cdot b\defeq \varphi(a)\cdot b\) restricts to a compact operator
\(f(a)\) on~\(\Hilm\) and \(f\colon A\to\Comp(\Hilm)\) is essential.  It
can be checked that this construction is functorial.

In the following proposition, we require that the category of
\(G\)\nbd{}\Cstar{}algebras \(\Cstarnice\) be closed under
\textsc{Morita}--\textsc{Rieffel} equivalence and consist of
\(\sigma\)\nbd{}unital \(G\)\nbd{}\Cstar{}algebras.  We let
\(\Corr{\Cstarnice}\) be the full subcategory of \(\Corr{G}\) with
object class~\(\Cstarnice\).

\begin{proposition}
  \label{pro:universal_stable_functor}
  The functor \(\natural\colon \Cstarnice \to \Corr{\Cstarnice}\) is the
  universal \Cstar{}stable functor on~\(\Cstarnice\); that is, it is
  \Cstar{}stable, and any other such functor factors uniquely
  through~\(\natural\).
\end{proposition}

\begin{proof}
  First we sketch the proof in the non-equivariant case.
  We verify that~\(\natural\) is \Cstar{}stable.  The
  \textsc{Morita}--\textsc{Rieffel} equivalence between
  \(\Comp(\ell^2\N)\otimes A\cong \Comp\bigl(\ell^2(\N,A)\bigr)\)
  and~\(A\) is implemented by the \textsc{Hilbert} module
  \(\ell^2(\N,A)\), which yields a correspondence
  \(\bigl(\ID,\ell^2(\N,A)\bigr)\) from \(\Comp(\ell^2\N)\otimes A\)
  to~\(A\); this is inverse to the correspondence induced by a corner
  embedding \(A\to \Comp(\ell^2\N)\otimes A\).

  A \textsc{Hilbert} \(B\)\nbd{}module~\(\Hilm\) with an essential
  \Star{}homomorphism \(A\to\Comp(\Hilm)\) is countably generated
  because~\(A\) is assumed \(\sigma\)\nbd{}unital.  \textsc{Kasparov}'s
  Stabilisation Theorem yields an isometric embedding \(\Hilm\to
  \ell^2(\N,B)\).  Hence we get \Star{}homomorphisms
  \[
  A\to \Comp(\ell^2\N)\otimes B \leftarrow B.
  \]
  This diagram induces a map \(F(A)\to F(\Comp(\ell^2\N)\otimes B)\cong
  F(B)\) for any \Cstar{}stable functor~\(F\).  Now we should check that
  this well-defines a functor \(\bar{F}\colon \Corr{\Cstarnice}\to
  \Cat\) with \(\bar{F}\circ\natural = F\), and that this yields the
  only such functor.  We omit these computations.

  The generalisation to the equivariant case uses the crucial property
  of the left regular representation that \(L^2(G)\otimes\Hils \cong
  L^2(G\times\N)\) for any countably infinite-dimensional
  \(G\)\nbd{}\textsc{Hilbert} space~\(\Hils\).  Since we replace \(A\)
  and~\(B\) by \(A_\Comp\) and~\(B_\Comp\) in the definition of
  correspondence right away, we can use this to repair a possible lack
  of \(G\)\nbd{}equivariance; similar ideas appear in~\cite{Meyer:KKG}.
\end{proof}

\begin{example}
  \label{exa:inner_auto}
  Let~\(u\) be a \(G\)\nbd{}invariant multiplier of~\(B\) with
  \(u^*u=1\); such~\(u\) are also called \emph{isometries}.  Then
  \(b\mapsto ubu^*\) defines a \Star{}homomorphism \(B\to B\).  The
  resulting correspondence \(B\dasharrow B\) is isomorphic as a
  correspondence to the identity correspondence: the isomorphism is
  given by left multiplication with~\(u\), which defines a
  \(G\)\nbd{}equivariant unitary operator from~\(B\) to the closure of
  \(uBu^*\cdot B = u\cdot B\).

  Hence inner endomorphisms act trivially on \Cstar{}stable functors.
  Actually, this is one of the computations that we have omitted in the
  proof above; the argument can be found
  in~\cite{Cuntz-Meyer-Rosenberg:MFO}.
\end{example}

Now we make the definition of a correspondence more concrete if \(A\) is
unital.  We have an essential \Star{}homomorphism \(\varphi\colon
A\to\Comp(\Hilm)\) for some \(G\)\nbd{}equivariant \textsc{Hilbert}
\(B\)\nbd{}module~\(\Hilm\).  Since~\(A\) is unital, this means that
\(\Comp(\Hilm)\) is unital and~\(\varphi\) is a unital
\Star{}homomorphism.  Then~\(\Hilm\) is finitely generated.  Thus
\(\Hilm= B^\infty\cdot p\) for some projection \(p\in\Mat_\infty(B)\)
and~\(\varphi\) is a \Star{}homomorphism \(\varphi\colon A\to
\Mat_\infty(B)\) with \(\varphi(1)=p\).  Two \Star{}homomorphisms
\(\varphi_1,\varphi_2\colon A\rightrightarrows \Mat_\infty(B)\) yield
isomorphic correspondences if and only if there is a partial isometry
\(v\in \Mat_\infty(B)\) with \(v\varphi_1(x)v^*=\varphi_2(x)\) and
\(v^*\varphi_2(a)v=\varphi_1(a)\) for all \(a\in A\).

Finally, we combine homotopy invariance and \Cstar{}stability and
consider the universal \Cstar{}stable homotopy-invariant functor.  This
functor is much easier to characterise: the morphisms in the resulting
universal category are simply the homotopy classes of
\(G\)\nbd{}equivariant \Star{}homomorphisms
\(\Comp\bigl(L^2(G\times\N)\bigr)\otimes A \to
\Comp\bigl(L^2(G\times\N)\bigr)\otimes B)\) (see
\cite{Meyer:KKG}*{Proposition 6.1}).  Alternatively, we get the same
category if we use homotopy classes of correspondences \(A\dasharrow B\)
instead.

\subsection{Exactness properties}
\label{sec:exactness}

Throughout this subsection, we consider functors \(F\colon
\Cstarnice\to\Cat\) with values in an exact category~\(\Cat\).
If~\(\Cat\) is merely additive to begin with, we can equip it with the
trivial exact category structure for which all extensions split.  We
also suppose that~\(\Cstarnice\) is closed under the kinds of
\Cstar{}algebra extensions that we consider; depending on the notion of
exactness, this means: direct product extensions, split extensions,
cp-split extensions, or all extensions, respectively.  Recall that split
extensions in~\(\CstarcatG\) are required to split by a
\(G\)\nbd{}equivariant \Star{}homomorphism.

\subsubsection{Additive functors}
\label{sec:additive_functor}

The most trivial split extensions in~\(\CstarcatG\) are the
\emph{product extensions} \(A\into A\times B\prto B\) for two objects
\(A,B\).  In this case, the coordinate embeddings and projections
provide maps
\begin{equation}
  \label{eq:direct_product_diagram}
  A\leftrightarrows A\times B \leftrightarrows B.
\end{equation}

\begin{definition}
  \label{def:additive}
  We call~\(F\) \emph{additive} if it maps product
  diagrams~\eqref{eq:direct_product_diagram} in~\(\Cstarnice\) to direct
  sum diagrams in~\(\Cat\).
\end{definition}

There is a partially defined addition on \Star{}homomorphisms: call two
parallel \Star{}homomorphisms \(\varphi,\psi\colon A\rightrightarrows
B\) \emph{orthogonal} if \(\varphi(a_1)\cdot \psi(a_2)=0\) for all
\(a_1,a_2\in A\).  Equivalently, \(\varphi+\psi\colon a\mapsto
\varphi(a)+\psi(a)\) is again a \Star{}homomorphism.

\begin{lemma}
  \label{lem:additive}
  The functor~\(F\) is additive if and only if, for all
  \(A,B\inOb\Cstarnice\), the maps \(\Hom(A,B)\to
  \Cat\bigl(F(A),F(B)\bigr)\) satisfy
  \(F(\varphi+\psi)=F(\varphi)+F(\psi)\) for all pairs of orthogonal
  parallel \Star{}homomorphisms \(\varphi,\psi\).
\end{lemma}

Alternatively, we may also require additivity for coproducts (that is,
free products).  Of course, this only makes sense if~\(\Cstarnice\) is
closed under coproducts in~\(\CstarcatG\).  The coproduct \(A\sqcup B\)
in~\(\CstarcatG\) comes with canonical maps \(A \leftrightarrows A\sqcup
B \leftrightarrows B\) as well; the maps \(\iota_A\colon A\to A\sqcup
B\) and \(\iota_B\colon B\to A\sqcup B\) are the coordinate embeddings,
the maps \(\pi_A\colon A\sqcup B\to A\) and \(\pi_B\colon A\sqcup B\to
B\) restrict to \((\ID_A,0)\) and \((0,\ID_B)\) on \(A\) and~\(B\),
respectively.

\begin{definition}
  \label{def:additive_co}
  We call~\(F\) \emph{additive on coproducts} if it maps coproduct
  diagrams \(A \leftrightarrows A\sqcup B \leftrightarrows B\) to direct
  sum diagrams in~\(\Cat\).
\end{definition}

The coproduct and product are related by a canonical
\(G\)\nbd{}equivariant \Star{}homomorphism \(\varphi\colon A\sqcup
B\prto A\times B\) that is compatible with the maps to and from \(A\)
and~\(B\), that is, \(\varphi\circ\iota_A=\iota_A\), \(\pi_A\circ
\varphi =\pi_A\), and similarly for~\(B\).  There is no map backwards,
but there is a \emph{correspondence} \(\psi\colon A\times B\dasharrow
A\sqcup B\), which is induced by the \(G\)\nbd{}equivariant
\Star{}homomorphism
\[
A\times B\to \Mat_2(A\sqcup B),\qquad
(a,b) \mapsto \begin{pmatrix}\iota_A(a)&0\\0&\iota_B(b)\end{pmatrix}.
\]
It is easy to see that the composite correspondence
\(\varphi\circ\psi\) is equal to the identity correspondence on
\(A\times B\).  The other composite \(\psi\circ\varphi\) is not the
identity correspondence, but it is homotopic to it (see
\cites{Cuntz:GenHom,Cuntz:New_Look}).  This yields:

\begin{proposition}
  \label{pro:additivity_compare}
  If~\(F\) is \Cstar{}stable and homotopy invariant, then the canonical
  map \(F(\varphi)\colon F(A\sqcup B) \to F(A\times B)\) is invertible.
  Therefore, additivity and additivity for coproducts are equivalent for
  such functors.
\end{proposition}

The correspondence~\(\psi\) exists because the stabilisation creates
enough room to replace \(\iota_A\) and~\(\iota_B\) by homotopic
homomorphisms with orthogonal ranges.  We can achieve the same effect by
a suspension (shift \(\iota_A\) and~\(\iota_B\) to the open intervals
\((0,\nicefrac{1}{2})\) and \((\nicefrac{1}{2},1)\), respectively).
Therefore, any homotopy invariant functor satisfies
\(F\bigl(\Sus(A\sqcup B)\bigr) \cong F\bigl(\Sus(A\times B)\bigr)\).

\subsubsection{Split-exact functors}
\label{sec:split-exact}

\begin{definition}
  \label{def:split-exact}
  We call~\(F\) \emph{split-exact} if, for any split extension
  \(K\overset{i}\into E \overset{p}\prto Q\) with section \(s\colon Q\to
  E\), the map \(\bigl(F(i),F(s)\bigr)\colon F(K)\oplus F(Q)\to F(E)\)
  is invertible.
\end{definition}

It is clear that split-exact functors are additive.

Split-exactness is useful because of the following construction of
\textsc{Joachim Cuntz} \cite{Cuntz:GenHom}.

Let \(B\triangleleft E\) be a \(G\)\nbd{}invariant ideal and let
\(f_+,f_-\colon A\rightrightarrows E\) be \(G\)\nbd{}equivariant
\Star{}homomorphisms with \(f_+(a)-f_-(a)\in B\) for all \(a\in A\).
Equivalently, \(f_+\) and~\(f_-\) both lift the same morphism
\(\bar{f}\colon A\to E/B\).  The data \((A,f_+,f_-,E,B)\) is called a
\emph{quasi-homomorphism} from~\(A\) to~\(B\).

Pulling back the extension \(B\into E\prto E/B\) along~\(\bar{f}\), we
get an extension \(B\into E'\prto A\) with two sections
\(f_+',f_-'\colon A\rightrightarrows E'\).  The split-exactness
of~\(F\) shows that \(F(B) \into F(E') \prto F(A)\) is a split
extension in~\(\Cat\).  Since both \(F(f_-')\) and \(F(f_+')\) are
sections for it, we get a map \(F(f_+')-F(f_-')\colon F(A)\to F(B)\).
\emph{Thus a quasi-homomorphism induces a map \(F(A)\to F(B)\)
  if~\(F\) is split-exact.}  The formal properties of this
construction are summarised in~\cite{Cuntz-Meyer-Rosenberg:MFO}.

Given a \Cstar{}algebra~\(A\), there is a \emph{universal
  quasi-homomorphism} out of~\(A\).  Let \(Q(A)\defeq A\sqcup A\) be the
coproduct of two copies of~\(A\) and let \(\pi_A\colon Q(A)\to A\) be
the \emph{folding homomorphism} that restricts to \(\ID_A\) on both
factors.  Let \(q(A)\) be its kernel.  The two canonical embeddings
\(A\to A\sqcup A\) are sections for the folding homomorphism.  Hence we
get a quasi-homomorphism \(A\rightrightarrows Q(A)\triangleright q(A)\).
The universal property of the free product shows that any
quasi-homomorphism yields a \(G\)\nbd{}equivariant \Star{}homomorphism
\(q(A)\to B\).

\begin{theorem}
  \label{the:stable_spexact_homotopy}
  Suppose~\(\Cstarnice\) is closed under split extensions and tensor
  products with \(\CONT([0,1])\) and \(\Comp(\ell^2\N)\).  If \(F\colon
  \Cstarnice\to\Cat\) is \Cstar{}stable and split-exact, then~\(F\) is
  homotopy invariant.
\end{theorem}

This is a deep result of \textsc{Nigel Higson}
\cite{Higson:Algebraic_K_stable}; a simple proof can be found in
\cite{Cuntz-Meyer-Rosenberg:MFO}.  Besides basic properties of
quasi-homomorphisms, it uses that inner endomorphisms act identically on
\Cstar{}stable functors (Example~\ref{exa:inner_auto}).

Actually, the literature only contains
Theorem~\ref{the:stable_spexact_homotopy} for functors on~\(\Cstarcat\).
But the proof in~\cite{Cuntz-Meyer-Rosenberg:MFO} works for functors on
categories~\(\Cstarnice\) as above.

\subsubsection{Exact functors}
\label{exact_functors}

\begin{definition}
  \label{def:exact}
  We call~\(F\) \emph{exact} if \(F(K)\to F(E)\to F(Q)\) is exact (at
  \(F(E)\)) for any extension \(K\into E\prto Q\) in~\(\Cstarnice\).
  More generally, given a class~\(\mathcal{E}\) of extensions
  in~\(\Cstarnice\) like, say, the class of equivariantly cp-split
  extensions, we define exactness for extensions in~\(\mathcal{E}\).
\end{definition}

It is easy to see that exact functors are additive.

Most functors we are interested in satisfy homotopy invariance and
\textsc{Bott} periodicity, and these two properties prevent a non-zero
functor from being exact in the stronger sense of being \emph{left} or
\emph{right} exact.  This explains why our notion of exactness is much
weaker than usual in homological algebra.

It is reasonable to require that a functor be part of a \emph{homology}
theory, that is, a sequence of functors \((F_n)_{n\in\Z}\) together with
natural long exact sequences for all extensions \cite{Schochet:Top3}.
We do not require this additional information because it tends to be
hard to get \emph{a priori} but often comes for free \emph{a
  posteriori}:

\begin{proposition}
  \label{pro:exact_longsequence}
  Suppose that~\(F\) is homotopy invariant and exact (or exact for
  equivariantly cp-split extensions).  Then~\(F\) has long exact
  sequences of the form
  \[
  \dotsb \to F\bigl(\Sus(K)\bigr) \to F\bigl(\Sus(E)\bigr) \to
  F\bigl(\Sus(Q)\bigr) \to F(K) \to F(E) \to F(Q)
  \]
  for any (equivariantly cp-split) extension \(K\into E\prto Q\).  In
  particular, \(F\) is split-exact.
\end{proposition}

See \textsection21.4 in~\cite{Blackadar:Book} for the proof.

There probably exist exact functors that are not split-exact.  It is
likely that the algebraic \(\K_1\)-functor provides a counterexample: it
exact but not split-exact on the category of rings
\cite{Rosenberg:Algebraic_K}; but I do not know a counterexample to its
split-exactness involving only \Cstar{}algebras.

Proposition~\ref{pro:exact_longsequence} and \textsc{Bott} periodicity
yield long exact sequences that are infinite in \emph{both} directions.
Thus an exact homotopy invariant functor that satisfies \textsc{Bott}
periodicity is part of a homology theory in a canonical way.

For universal constructions, we should replace a single functor by a
\emph{homology theory}, that is, a sequence of functors.  The
universal functors in this context are non-stable versions of
\(\E\)\nbd{}theory and \(\KK\)-theory.  We refer to
\cite{Houghton-Larsen-Thomsen} for details.

A weaker property than exactness is the existence of
\emph{\textsc{Puppe} exact sequences} for mapping cones.  The
\textsc{Puppe} exact sequence is the special case of the long exact
sequence of Proposition~\ref{pro:exact_longsequence} for
extensions of the form \(\Sus(B) \into \Cone(f)\prto A\) for a morphism
\(f\colon A\to B\).  In practice, the exactness of a functor is often
established by reducing it to the \textsc{Puppe} exact sequence.  Let
\(K\overset{i}\into E\overset{p}\prto Q\) be an extension.  A variant of
the \textsc{Puppe} exact sequence yields the long exact sequence for the
extension \(\Cone(p)\into \Cyl(p)\prto Q\).  There is a canonical
morphism of extensions
\[\xymatrix{
  K\ \ar@{>->}[r]\ar[d]&
  E \ar@{->>}[r]\ar[d]&
  Q\ar@{=}[d]\\
  \Cone(p)\ \ar@{>->}[r]&
  \Cyl(p) \ar@{->>}[r]&
  Q,
}
\]
where the vertical map \(E\to\Cyl(p)\) is a homotopy equivalence.  Hence
a functor with \textsc{Puppe} exact sequences is exact for \(K\into
E\prto Q\) if and only if it maps the vertical map \(K\to\Cone(p)\) to
an isomorphism.

\section{\textsc{Kasparov} theory}
\label{sec:KK_E}

We define \(\KK^G\) as the universal split-exact \Cstar{}stable functor
on \(\CstarsepG\); since split-exact and \Cstar{}stable functors are
automatically homotopy invariant, \(\KK^G\) is the universal split-exact
\Cstar{}stable homotopy functor as well.  The universal property of
\textsc{Kasparov} theory due to \textsc{Higson} and \textsc{Cuntz}
asserts that this is equivalent to \textsc{Kasparov}'s definition.  We
examine some basic properties of \textsc{Kasparov} theory and, in
particular, show how to get functors between \textsc{Kasparov}
categories.

We let~\(\E^G\) be the universal exact \Cstar{}stable homotopy functor
on \(\CstarsepG\) or, equivalently, the universal exact, split-exact,
and \Cstar{}stable functor.

\textsc{Kasparov}'s own definition of his theory is inspired by previous
work of \textsc{Atiyah} \cite{Atiyah:Ell} on \(\K\)\nbd{}homology;
later, he also interacted with the work of
\textsc{Brown}--\textsc{Douglas}--\textsc{Fillmore} \cite{BDF} on
extensions of \Cstar{}algebras.  A construction in abstract homotopy
theory provides a homology theory for spaces that is dual to
\(\K\)\nbd{}theory.  \textsc{Atiyah} realized that certain abstract
elliptic differential operators provide cycles for this dual theory; but
he did not know the equivalence relation to put on these cycles.
\textsc{Brown}--\textsc{Douglas}--\textsc{Fillmore} studied extensions
of \(\CONT_0(X)\) (and more general \Cstar{}algebras) by the compact
operators and found that the resulting structure set is naturally
isomorphic to a \(\K\)\nbd{}homology group.

\textsc{Kasparov} unified and vastly generalised these two results,
defining a bivariant functor \(\KK_*(A,B)\) that combines
\(\K\)\nbd{}theory and \(\K\)\nbd{}homology and that is closely related
to the classification of extensions \(B\otimes\Comp\into E\prto A\) (see
\cites{Kasparov:Elliptic, Kasparov:KK}).  A deep theorem of
\textsc{Kasparov} shows that two reasonable equivalence relations for
these cycles coincide; this clarifies the homotopy invariance of the
extension groups of \textsc{Brown}--\textsc{Douglas}--\textsc{Fillmore}.
Furthermore, he constructed an equivariant version of his theory in
\cite{Kasparov:Novikov} and applied it to prove the \textsc{Novikov}
conjecture for discrete subgroups of \textsc{Lie} groups.

The most remarkable feature of \textsc{Kasparov} theory is an
associative product on \(\KK\) called \emph{\textsc{Kasparov} product}.
This generalises various known product constructions in
\(\K\)\nbd{}theory and \(\K\)\nbd{}homology and allows to view \(\KK\)
as a category.

In applications, we usually need some non-obvious \(\KK\)-element, and
we must compute certain \textsc{Kasparov} products explicitly.  This
requires a concrete description of \textsc{Kasparov} cycles and their
products.  Since both are somewhat technical, we do not discuss them
here and merely refer to~\cite{Blackadar:Book} for a detailed
treatment and to~\cite{Skandalis:Survey_KK} for a very useful survey
article.  Instead, we use \textsc{Higson}'s characterisation of
\(\KK\) by a universal property \cite{Higson:Characterize_KK}, which
is based on ideas of \textsc{Cuntz} (\cites{Cuntz:New_Look,
  Cuntz:GenHom}).  The extension to the equivariant case is due to
\textsc{Thomsen} \cite{Thomsen:KKG_universal}.  A simpler proof of
\textsc{Thomsen}'s theorem and various related results can be found
in~\cite{Meyer:KKG}.

We do not discuss \(\KK^G\) for \(\Ztwo\)-graded
\(G\)\nbd{}\Cstar{}algebras here because it does not fit so well with
the universal property approach, which would simply yield
\(\KK^{G\times\Ztwo}\) because \(\Ztwo\)-graded
\(G\)\nbd{}\Cstar{}algebras are the same as
\(G\times\Ztwo\)-\Cstar{}algebras.  The relationship between the two
theories is explained in~\cite{Meyer:KKG}, following \textsc{Ulrich
  Haag} \cite{Haag:Graded_KK2}.  The graded version of \textsc{Kasparov}
theory is often useful because it allows us to treat even and odd
\(\KK\)-cycles simultaneously.

Fix a locally compact group~\(G\).  The \textsc{Kasparov} groups
\(\KK_0^G(A,B)\) for \(A,B\in\CstarsepG\) form the morphisms sets \(A\to
B\) of a category, which we denote by \(\KKcat^G\); the composition in
\(\KKcat^G\) is the \textsc{Kasparov} product.  The categories
\(\CstarsepG\) and \(\KKcat^G\) have the same objects.  We have a
canonical functor
\[
\KK^G\colon \CstarsepG\to\KKcat^G
\]
that acts identically on objects.  This functor contains all information
about equivariant \textsc{Kasparov} theory for~\(G\).

\begin{definition}
  \label{def:KK_equivalence}
  A \(G\)\nbd{}equivariant \Star{}homomorphism \(f\colon A\to B\) is
  called a \emph{\(\KK^G\)-equivalence} if \(\KK^G(f)\) is invertible in
  \(\KKcat^G\).
\end{definition}

\begin{theorem}
  \label{the:KK_universal}
  Let \(\Cstarnice\) be a full subcategory of \(\CstarsepG\) that is
  closed under suspensions, \(G\)\nbd{}equivariantly cp-split
  extensions, and \textsc{Morita}--\textsc{Rieffel} equivalence.  Let
  \(\KKcat^G(\Cstarnice)\) be the full subcategory of \(\KKcat^G\) with
  object class~\(\Cstarnice\) and let
  \(\KK^G_\Cstarnice\colon \Cstarnice\to\KKcat^G(\Cstarnice)\) be the
  restriction of \(\KK^G\).

  The functor \(\KK^G_\Cstarnice\colon
  \Cstarnice\to\KKcat^G_\Cstarnice\) is the universal split-exact
  \Cstar{}stable functor; in particular, \(\KKcat^G(\Cstarnice)\) is an
  additive category.  In addition, it has the following properties and
  is, therefore, universal among functors on~\(\Cstarnice\) with some of
  these extra properties: it is
  \begin{itemize}
  \item homotopy invariant;

  \item exact for \(G\)\nbd{}equivariantly cp-split extensions;

  \item satisfies \textsc{Bott} periodicity, that is, in \(\KKcat^G\)
    there are natural isomorphisms \(\Sus^2(A) \cong A\) for all
    \(A\inOb\KKcat^G\).

  \end{itemize}
\end{theorem}

\begin{corollary}
  \label{cor:Bott}
  Let \(F\colon \Cstarnice\to\Cat\) be split-exact and \Cstar{}stable.
  Then~\(F\) factors uniquely through \(\KK^G_\Cstarnice\), is homotopy
  invariant, and satisfies \textsc{Bott} periodicity.  A
  \(\KK^G\)-equivalence \(A\to B\) in \(\Cstarnice\) induces an
  isomorphism \(F(A)\to F(B)\).
\end{corollary}

We will view the universal property of Theorem~\ref{the:KK_universal} as
a definition of \(\KKcat^G\) and thus of the groups \(\KK_0^G(A,B)\).
We also let
\[
\KK_n^G(A,B) \defeq \KK^G\bigl(A,\Sus^n(B)\bigr);
\]
since the \textsc{Bott} periodicity isomorphism identifies
\(\KK_2^G\cong \KK_0^G\), this yields a \(\Ztwo\)\nbd{}graded theory.

Now we describe \(\KK^G_0(A,B)\) more concretely.  Recall
\(A_\Comp\defeq A\otimes \Comp(L^2G)\).

\begin{proposition}
  \label{pro:KKconcrete}
  Let \(A\) and~\(B\) be two \(G\)\nbd{}\Cstar{}algebras.  There is a
  natural bijection between the morphism sets \(\KK_0^G(A,B)\) in
  \(\KKcat^G\) and the set \([q(A_\Comp), B_\Comp\otimes
  \Comp(\ell^2\N)]\) of homotopy classes of \(G\)\nbd{}equivariant
  \Star{}homomorphisms from \(q(A_\Comp)\) to \(B_\Comp\otimes
  \Comp(\ell^2\N)\).
\end{proposition}

\begin{proof}
  The canonical functor \(\CstarsepG\to\KKcat^G\) is \Cstar{}stable and
  split-exact, and therefore homotopy invariant by Theorem
  \ref{the:stable_spexact_homotopy} (this is already asserted in
  Theorem~\ref{the:KK_universal}).  Proposition
  \ref{pro:additivity_compare} yields that it is additive for
  coproducts.  Split-exactness for the split extension \(q(A)\into
  Q(A)\prto A\) shows that \(\ID_A*0\colon Q(A)\to A\) restricts to a
  \(\KK^G\)-equivalence \(q(A)\sim A\).  Similarly, \Cstar{}stability
  yields \(\KK^G\)-equivalences \(A\sim A_\Comp\) and \(B\sim
  B_\Comp\otimes \Comp(\ell^2\N)\).  Hence homotopy classes of
  \Star{}homomorphisms from \(q(A_\Comp)\) to \(B_\Comp\otimes
  \Comp(\ell^2\N)\) yield classes in \(\KK_0^G(A,B)\).  Using the
  concrete description of \textsc{Kasparov} cycles, which we have not
  discussed, it is checked in~\cite{Meyer:KKG} that this map yields a
  bijection as asserted.
\end{proof}

Another equivalent description is
\[
\KK_0^G(A,B) \cong
[q(A_\Comp)\otimes \Comp(\ell^2\N), q(B_\Comp)\otimes \Comp(\ell^2\N)];
\]
in this approach, the \textsc{Kasparov} product becomes simply the
composition of morphisms.  Proposition \ref{pro:KKconcrete} suggests
that \(q(A_\Comp)\) and \(B_\Comp\otimes \Comp(\ell^2\N)\) may be the
cofibrant and fibrant replacement of \(A\) and~\(B\) in some model
category related to \(\KK^G\).  But it is not clear whether this is the
case.  The model category structure constructed
in~\cite{Joachim-Johnson:Model} is certainly quite different.

By the universal property, \(\K\)\nbd{}theory descends to a functor on
\(\KKcat\), that is, we get canonical maps
\[
\KK_0(A,B) \to \Hom\bigl(\K_*(A),\K_*(B)\bigr)
\]
for all separable \Cstar{}algebras \(A,B\), where the right hand side
denotes grading-preserving group homomorphisms.  For \(A=\C\), this
yields a map \(\KK_0(\C,B)\to \Hom\bigl(\Z,\K_0(B)\bigr) \cong
\K_0(B)\).  Using suspensions, we also get a corresponding map
\(\KK_1(\C,B)\to \K_1(B)\).

\begin{theorem}
  \label{the:KK_versus_K}
  The maps \(\KK_*(\C,B)\to \K_*(B)\) constructed above are isomorphisms
  for all \(B\inOb\Cstarsep\).
\end{theorem}

Thus \textsc{Kasparov} theory is a bivariant generalisation of
\(\K\)\nbd{}theory.  Roughly speaking, \(\KK_*(A,B)\) is the place where
maps between \(\K\)\nbd{}theory groups live.  Most constructions of such
maps, say, in index theory can, in fact, be improved to yield elements
of \(\KK_*(A,B)\).  One reason why this has to be so is the
\emph{Universal Coefficient Theorem} (UCT), which computes
\(\KK_*(A,B)\) from \(\K_*(A)\) and \(\K_*(B)\) for many
\Cstar{}algebras \(A,B\).  If~\(A\) satisfies the UCT, then any grading
preserving group homomorphism \(\K_*(A)\to\K_*(B)\) lifts to an element
of \(\KK_0(A,B)\).

\subsection{Extending functors and identities to
  \texorpdfstring{$\KKcat^G$}{Kasparov theory}}
\label{sec:construct_functors_universal}

We can use the universal property to extend various functors
\(\CstarsepG\to\CstarsepG[H]\) to functors \(\KKcat^G\to\KKcat^H\).  We
explain this by an example:

\begin{proposition}
  \label{pro:construct_descent}
  The full and reduced crossed product functors
  \[
  G\rcross\blank, G\cross\blank\colon \CstarcatG\to\Cstarcat
  \]
  extend to functors \(G\rcross\blank, G\cross\blank\colon
  \KKcat^G\to\KKcat\) called \emph{descent functors}.
\end{proposition}

\textsc{Gennadi Kasparov} \cite{Kasparov:Novikov} constructs these
functors directly using the concrete description of \textsc{Kasparov}
cycles.  This requires a certain amount of work; in particular, checking
functoriality involves knowing how to compute \textsc{Kasparov}
products.  The construction via the universal property is formal:

\begin{proof}
  We only write down the argument for \emph{reduced} crossed products,
  the other case is similar.  It is well-known that \(G\rcross
  \bigl(A\otimes \Comp(\Hils)\bigr) \cong (G\rcross A)\otimes
  \Comp(\Hils)\) for any \(G\)\nbd{}\textsc{Hilbert} space~\(\Hils\).
  Therefore, the composite functor
  \[
  \CstarsepG\xrightarrow{G\rcross\blank}
  \Cstarsep \xrightarrow{\KK}
  \KKcat
  \]
  is \Cstar{}stable.  Proposition \ref{pro:exact_on_cpsplit} shows that
  this functor is split-exact as well (regardless of whether~\(G\) is an
  exact group).  Now the universal property provides an extension to a
  functor \(\KKcat^G\to \KKcat\).
\end{proof}

Similarly, we get functors
\[
A\minotimes\blank, A\maxotimes\blank\colon \KKcat^G\to\KKcat^G
\]
for any \(G\)\nbd{}\Cstar{}algebra~\(A\).  Since these extensions are
natural, we even get bifunctors
\[
\minotimes,\maxotimes\colon \KKcat^G \times \KKcat^G\to\KKcat^G.
\]
The associativity, commutativity, and unit constraints in \(\CstarcatG\)
induce corresponding constraints in \(\KKcat^G\), so that both
\(\minotimes\) and~\(\maxotimes\) turn \(\KKcat^G\) into a symmetric
monoidal category.

Another example is the functor \(\tau\colon \Cstarcat\to\CstarcatG\)
that equips a \Cstar{}algebra with the trivial \(G\)\nbd{}action; it
extends to a functor \(\tau\colon \KKcat\to\KKcat^G\).

The universal property also allows us to \emph{prove identities} between
functors.  For instance, we have natural isomorphisms \(G\rcross
(\tau(A)\minotimes B) = A\minotimes (G\rcross B)\) for all
\(G\)\nbd{}\Cstar{}algebras~\(B\).  To begin with, naturality means that
the diagram
\[\xymatrix{
  G\rcross (\tau(A_1)\minotimes B_1) \ar[r]^{\cong}
  \ar[d]_{G\rcross \tau(\alpha)\minotimes \beta} &
  A_1\minotimes (G\rcross B_1)
  \ar[d]^{\alpha\minotimes (G\rcross \beta)} \\
  G\rcross (\tau(A_2)\minotimes B_2) \ar[r]^{\cong} &
  A_2\minotimes (G\rcross B_2) \\
}
\]
commutes if \(\alpha\colon A_1\to A_2\) and \(\beta\colon B_1\to B_2\)
are a \Star{}homomorphism and a \(G\)\nbd{}equivariant
\Star{}homomorphism, respectively.  Two applications of the uniqueness
part of the universal property show that this diagram remains
commutative in \(\KKcat\) if \(\alpha\in\KK_0(A_1,A_2)\) and
\(\beta\in\KK_0^G(B_1,B_2)\).  Similar remarks apply to the natural
isomorphism \(G\cross(\tau(A)\maxotimes B) \cong A\maxotimes (G\cross
B)\) and hence to the isomorphisms \(G\cross \tau(A) \cong \Cfull(G)
\maxotimes A\) and \(G\rcross \tau(A) \cong \Cred(G) \minotimes A\).

\emph{Adjointness} relations in \textsc{Kasparov} theory are usually
proved most easily by constructing the unit and counit of the
adjunction.  For instance, if~\(G\) is a compact group then the
functor~\(\tau\) is left adjoint to \(G\cross\blank= G\rcross\blank\),
that is, for all \(A\inOb\KKcat\) and \(B\inOb\KKcat^G\), we have
natural isomorphisms
\begin{equation}
  \label{eq:tau_cross_adjoint}
  \KK^G_*(\tau(A),B) \cong \KK_*(A,G\cross B).
\end{equation}
This is also known as the \emph{\textsc{Green}--\textsc{Julg} Theorem}.
For \(A=\C\), it specialises to a natural isomorphism \(\K^G_*(B) \cong
\K_*(G\cross B)\); this was one of the first appearances of
non-commutative algebras in topological \(\K\)\nbd{}theory.

\begin{proof}[Proof of~\eqref{eq:tau_cross_adjoint}]
  We already know that \(\tau\) and \(G\cross\blank\) are functors
  between \(\KKcat\) and \(\KKcat^G\).  It remains to construct natural
  elements
  \[
  \alpha_A\in\KK_0\bigl(A,G\cross \tau(A)\bigr), \qquad
  \beta_B\in\KK^G_0(\tau(G\cross B),B)
  \]
  for all \(A\inOb\KKcat\), \(B\inOb\KKcat^G\) that satisfy the
  conditions for unit and counit of adjunction
  \cite{MacLane:Categories}.

  The main point is that \(\tau(G\cross B)\) is the \(G\)\nbd{}fixed
  point subalgebra of \(B_\Comp = B\otimes\Comp(L^2G)\).  The embedding
  \(\tau(G\cross B) \to B_\Comp\) provides a \(G\)\nbd{}equivariant
  correspondence~\(\beta_B\) from \(\tau(G\cross B)\) to~\(B\) and thus
  an element of \(\KK^G_0(\tau(G\cross B),B)\).  This construction is
  certainly natural for \(G\)\nbd{}equivariant \Star{}homomorphisms and
  hence for \(\KK^G\)-morphisms by the uniqueness part of the universal
  property of \(\KK^G\).

  Let \(e_\tau\colon \C\to\Cfull(G)\) be the embedding that corresponds
  to the trivial representation of~\(G\).  Recall that \(G\cross
  \tau(A)\cong \Cfull(G)\otimes A\).  Hence the exterior product of the
  identity map on~\(A\) and \(\KK(e_\tau)\) provides \(\alpha_A\in
  \KK_0\bigl(A,G\cross \tau(A)\bigr)\).  Again, naturality for
  \Star{}homomorphisms is clear and implies naturality for morphisms in
  \(\KK\).

  Finally, it remains to check that
  \begin{gather*}
    \tau(A) \xrightarrow{\tau(\alpha_A)}
    \tau\bigl(G\cross\tau(A)\bigr) \xrightarrow{\beta_{\tau(A)}}
    \tau(A)\\
    G\cross B \xrightarrow{\alpha_{G\cross B}}
    G\cross \tau(G\cross B) \xrightarrow{G\cross\beta_B}
    G\cross B
  \end{gather*}
  are the identity morphisms in \(\KK^G\).  Then we get the desired
  adjointness using a general construction from category theory
  (see~\cite{MacLane:Categories}).  In fact, both composites are equal
  to the identity already as \emph{correspondences}, so that we do not
  have to know anything about \textsc{Kasparov} theory except its
  \Cstar{}stability to check this.
\end{proof}

A similar argument yields an adjointness relation
\begin{equation}
  \label{eq:cross_tau_adjoint}
  \KK^G_0\bigl(A,\tau(B)\bigr) \cong \KK_0(G\cross A,B)  
\end{equation}
for a \emph{discrete} group~\(G\).  More conceptually,
\eqref{eq:cross_tau_adjoint} corresponds via
\textsc{Baaj}--\textsc{Skandalis} duality \cite{Baaj-Skandalis:KKHopf}
to the \textsc{Green}--\textsc{Julg} Theorem for the dual quantum group
of~\(G\), which is \emph{compact} because~\(G\) is discrete.  But we can
also write down unit and counit of adjunction directly.

The trivial representation \(\Cfull(G)\to\C\) yields natural
\Star{}homomorphisms
\[
G\cross \tau(B) \cong \Cfull(G)\maxotimes B \to B
\]
and hence \(\beta_B\in \KK_0(G\cross\tau(B),B)\).  The canonical
embedding \(A\to G\cross A\) is \(G\)\nbd{}equivariant if we let~\(G\)
act on \(G\cross A\) by conjugation; but this action is inner, so that
\(G\cross A\) and \(\tau(G\cross A)\) are \(G\)\nbd{}equivariantly
\textsc{Morita}--\textsc{Rieffel} equivalent.  Thus the canonical
embedding \(A\to G\cross A\) yields a correspondence \(A\dasharrow
\tau(G\cross A)\) and \(\alpha_A\in \KK^G_0\bigl(A,\tau(G\cross
A)\bigr)\).  We must check that the composites
\begin{gather*}
  G\cross A \xrightarrow{G\cross\alpha_A}
  G\cross \tau(G\cross A) \xrightarrow{\beta_{G\cross A}}
  G\cross A,\\
  \tau(B) \xrightarrow{\alpha_{\tau(B)}}
  \tau\bigl(G\cross\tau(B)\bigr) \xrightarrow{\tau(\beta_B)}
  \tau(B)
\end{gather*}
are identity morphisms in \(\KK\) and \(\KK^G\), respectively.  Once
again, this holds already on the level of correspondences.

\subsection{Triangulated category structure}
\label{sec:tri_structure}

We turn \(\KKcat^G\) into a triangulated category by extending standard
constructions for topological spaces \cite{Meyer-Nest:BC}.  Some arrows
change direction because the functor~\(\CONT_0\) from spaces to
\Cstar{}algebras is contravariant.  We have already observed that
\(\KKcat^G\) is additive.  The suspension is given by
\(\Sigma^{-1}(A)\defeq \Sus(A)\).  Since \(\Sus^2(A)\cong A\) in
\(\KKcat^G\) by \textsc{Bott} periodicity, we have
\(\Sigma=\Sigma^{-1}\).  Thus we do not need formal desuspensions as for
the stable homotopy category.

\begin{definition}
  \label{def:triangle_exact}
  A triangle \(A\to B\to C\to \Sigma A\) in \(\KKcat^G\) is called
  \emph{exact} if it is isomorphic as a triangle to the \emph{mapping
    cone triangle}
  \[
  \Sus(B) \to \Cone(f) \to A \xrightarrow{f} B
  \]
  for some \(G\)\nbd{}equivariant \Star{}homomorphism~\(f\).
\end{definition}

Alternatively, we can use \(G\)\nbd{}equivariantly cp-split extensions
in \(\CstarsepG\).  Any such extension \(I\into E\prto Q\) determines a
class in \(\KK^G_1(Q,I)\cong \KK^G_0(\Sus(Q),I)\), so that we get a
triangle \(\Sus(Q)\to I\to E\to Q\) in \(\KKcat^G\).  Such triangles are
called \emph{extension triangles}.  A triangle in \(\KKcat^G\) is exact
if and only if it is isomorphic to the extension triangle of a
\(G\)\nbd{}equivariantly cp-split extension \cite{Meyer-Nest:BC}.

\begin{theorem}
  \label{the:KK_triangulated}
  With the suspension automorphism and exact triangles defined above,
  \(\KKcat^G\) is a triangulated category.  So is
  \(\KKcat^G(\Cstarnice)\) if \(\Cstarnice\subseteq\CstarsepG\) is
  closed under suspensions, \(G\)\nbd{}equivariantly cp-split
  extensions, and \textsc{Morita}--\textsc{Rieffel} equivalence as in
  Theorem~\ref{the:KK_universal}.
\end{theorem}

\begin{proof}
  That \(\KKcat^G\) is triangulated is proved in detail
  in~\cite{Meyer-Nest:BC}.  We do not discuss the triangulated category
  axioms here.  Most of them amount to properties of mapping cone
  triangles that can be checked by copying the corresponding arguments
  for the stable homotopy category (and reverting arrows).  These axioms
  hold for \(\KKcat^G(\Cstarnice)\) because they hold for \(\KKcat^G\).
  The only axiom that requires more care is the existence axiom for
  exact triangles; it requires any morphism to be part of an exact
  triangle.  We can prove this as in~\cite{Meyer-Nest:BC} using the
  concrete description of \(\KK^G_0(A,B)\) in
  Proposition~\ref{pro:KKconcrete}.  For some applications like the
  generalisation to \(\KKcat^G(\Cstarnice)\), it is better to use
  extension triangles instead.  Any \(f\in \KK^G_0(A,B)\cong
  \KK^G_1(\Sus(A),B)\) can be represented by a \(G\)\nbd{}equivariantly
  cp-split extension \(\Comp(\Hils_B)\into E\prto \Sus(A)\),
  where~\(\Hils_B\) is a full \(G\)\nbd{}equivariant \textsc{Hilbert}
  \(B\)\nbd{}module, so that \(\Comp(\Hils_B)\) is
  \(G\)\nbd{}equivariantly \textsc{Morita}--\textsc{Rieffel} equivalent
  to~\(B\).  The extension triangle of this extension contains~\(f\) and
  belongs to \(\KKcat^G(\Cstarnice)\) by our assumptions
  on~\(\Cstarnice\).
\end{proof}

Since model category structures related to \Cstar{}algebras are rather
hard to get (compare~\cite{Joachim-Johnson:Model}), triangulated
categories seem to provide the most promising formal setup for extending
results from classical spaces to \Cstar{}algebras.  An earlier attempt
can be found in~\cite{Schochet:Top3}.  Triangulated categories clarify
the basic bookkeeping with long exact sequences.
\textsc{Mayer}--\textsc{Vietoris} exact sequences and inductive limits
are discussed from this point of view in~\cite{Meyer-Nest:BC}.  More
importantly, this framework sheds light on more advanced constructions
like the \textsc{Baum}--\textsc{Connes} assembly map.  We will briefly
discuss this below.

\subsection{The Universal Coefficient Theorem}
\label{sec:UCT}

There is a very close relationship between \(\K\)\nbd{}theory and
\textsc{Kasparov} theory.  We have already seen that \(\K_*(A)\cong
\KK_*(\C,A)\) is a special case of \(\KK\).  Furthermore, \(\KK\)
inherits deep properties of \(\K\)\nbd{}theory such as \textsc{Bott}
periodicity.  Thus we may hope to express \(\KK_*(A,B)\) using only the
\(\K\)\nbd{}theory of \(A\) and~\(B\) \mdash{} at least for many \(A\)
and~\(B\).  This is the point of the \emph{Universal Coefficient
  Theorem}.

The \textsc{Kasparov} product provides a canonical homomorphism of
graded groups
\[
\gamma\colon \KK_*(A,B) \to \Hom_*\bigl(\K_*(A),\K_*(B)\bigr),
\]
where \(\Hom_*\) denotes the \(\Ztwo\)-graded \textsc{Abel}ian group of
all group homomorphisms \(\K_*(A)\to\K_*(B)\).  There are topological
reasons why~\(\gamma\) cannot always be invertible: since \(\Hom_*\) is
not exact, the bifunctor \(\Hom_*\bigl(\K_*(A),\K_*(B)\bigr)\) would not
be exact on cp-split extensions.  A construction of \textsc{Lawrence
  Brown} provides another natural map
\[
\kappa\colon \ker \gamma\to \Ext_*\bigl(\K_{*+1}(A),\K_*(B)\bigr).
\]
The following theorem is due to \textsc{Jonathan Rosenberg} and
\textsc{Claude Schochet} \cites{Schochet:Top2, Rosenberg-Schochet:UCT};
see also~\cite{Blackadar:Book}.

\begin{theorem}
  \label{the:bootstrap}
  The following are equivalent for a separable \Cstar{}algebra~\(A\):
  \begin{enumerate}[(a)]
  \item \(\KK_*(A,B)=0\) for all \(B\inOb\KKcat\) with \(\K_*(B)=0\);

  \item the map~\(\gamma\) is surjective and~\(\kappa\) is bijective for
    all \(B\inOb\KKcat\);

  \item for all \(B\inOb\KKcat\), there is a short exact sequence of
    \(\Ztwo\)\nbd{}graded \textsc{Abel}ian groups
    \[
    \Ext_*\bigl(\K_{*+1}(A),\K_*(B)\bigr) \into
    \KK_*(A,B) \prto
    \Hom_*\bigl(\K_*(A),\K_*(B)\bigr).
    \]

  \item \(A\) belongs to the smallest class of \Cstar{}algebras that
    contains~\(\C\) and is closed under \(\KK\)-equivalence,
    suspensions, countable direct sums, and cp-split extensions;

  \item \(A\) is \(\KK\)-equivalent to \(\CONT_0(X)\) for some pointed
    compact metrisable space~\(X\).

  \end{enumerate}
  If these conditions are satisfied, then the extension in~(c) is
  natural and splits, but the section is not natural.
\end{theorem}

The class of \Cstar{}algebras with these properties is also called the
\emph{bootstrap} class because of description~(d).  Alternatively, we
may say that they satisfy the \emph{Universal Coefficient Theorem}
because of~(c).  Since commutative \Cstar{}algebras are nuclear, (e)
implies that the natural map \(A\maxotimes B \to A\minotimes B\) is a
\(\KK\)-equivalence if \(A\) or~\(B\) belongs to the bootstrap class
\cite{Skandalis:K-nuclear}.  This fails for some~\(A\), so that the
Universal Coefficient Theorem does not hold for all~\(A\).  Remarkably,
this is the only obstruction to the Universal Coefficient Theorem known
at the moment: we know no nuclear \Cstar{}algebra that does not satisfy
the Universal Coefficient Theorem.  As a result, we can express
\(\KK_*(A,B)\) using only \(\K_*(A)\) and \(\K_*(B)\) for many \(A\)
and~\(B\).

When we restrict attention to nuclear \Cstar{}algebras, then the
bootstrap class is closed under various operations like tensor products,
arbitrary extensions and inductive limits (without requiring any
cp-sections), and under crossed products by torsion-free amenable
groups.  Remarkably, there are no general results about crossed products
by finite groups.

The Universal Coefficient Theorem and the universal property of \(\KK\)
imply that very few homology theories for (pointed compact metrisable)
spaces can extend to the non-commutative setting.  More precisely, if we
require the extension to be split-exact, \Cstar{}stable, and additive
for countable direct sums, then only \(\K\)\nbd{}theory with
coefficients is possible.  Thus we rule out most of the difficult (and
interesting) problems in stable homotopy theory.  But if we only want to
study \(\K\)\nbd{}theory, anyway, then the operator algebraic framework
usually provides very good analytical tools.  This is most valuable for
equivariant generalisations of \(\K\)\nbd{}theory.

\textsc{Jonathan Rosenberg} and \textsc{Claude Schochet}
\cite{Rosenberg-Schochet} have also constructed a spectral sequence
that, in favourable cases, computes \(\KK^G(A,B)\) from \(\K^G_*(A)\)
and \(\K^G_*(B)\); they require~\(G\) to be a compact \textsc{Lie} group
with torsion-free fundamental group and \(A\) and~\(B\) to belong to a
suitable bootstrap class.  This equivariant UCT is clarified in
\cites{Meyer-Nest:Homology_in_KK, Meyer-Nest:BC_Coaction}.

\subsection{E-Theory and asymptotic morphisms}
\label{sec:asymp_morph}

Recall that \textsc{Kasparov} theory is only exact for (equivariantly)
cp-split extensions.  \(\E\)\nbd{}Theory is a similar theory that is
exact for all extensions.

\begin{definition}
  \label{def:E}
  We let \(\E^G\colon \CstarsepG\to\Ecat^G\) be the universal
  \Cstar{}stable, exact homotopy functor.
\end{definition}

\begin{lemma}
  \label{lem:E_properties}
  The functor~\(\E^G\) is split-exact and factors through \(\KK^G\colon
  \CstarsepG\to\KK^G\).  Hence it satisfies \textsc{Bott} periodicity.
\end{lemma}

\begin{proof}
  Proposition~\ref{pro:exact_longsequence} shows that any exact homotopy
  functor is split-exact.  The remaining assertions now follow from
  Corollary~\ref{cor:Bott}.
\end{proof}

The functor~\(\E\) (for trivial~\(G\)) was first defined as above by
\textsc{Nigel Higson} \cite{Higson:Fractions}.  Then \textsc{Alain
  Connes} and \textsc{Nigel Higson} \cite{Connes-Higson:E} found a more
concrete description using \emph{asymptotic morphisms}.  This is what
made the theory usable.  The equivariant generalisation of the theory is
due to \textsc{Erik Guentner}, \textsc{Nigel Higson}, and \textsc{Jody
  Trout} \cite{Guentner-Higson-Trout}.

We write \(\E^G_n(A,B)\) for the space of morphisms \(A\to\Sus^n(B)\) in
\(\Ecat^G\).  \textsc{Bott} periodicity shows that there are only two
different groups to consider.

\begin{definition}
  \label{def:asymp_mor}
  The \emph{asymptotic algebra} of a \Cstar{}algebra~\(B\) is the
  \Cstar{}algebra
  \[
  \Asymp(B) \defeq \CONT_b(\R_+,B)/\CONT_0(\R_+,B).
  \]
  An \emph{asymptotic morphism} \(A\to B\) is a \Star{}homomorphism
  \(f\colon A\to\Asymp(B)\).
\end{definition}

Representing elements of \(\Asymp(B)\) by bounded functions
\([0,\infty)\to B\), we can represent~\(f\) by a family of maps
\(f_t\colon A\to B\) such that \(f_t(a)\in\CONT_b(\R_+,B)\) for each
\(a\in A\) and the map \(a\mapsto f_t(a)\) satisfies the conditions for
a \Star{}homomorphism asymptotically for \(t\to\infty\).  This provides
a concrete description of asymptotic morphisms and explains the name.

If a locally compact group~\(G\) acts on~\(B\), then \(\Asymp(B)\)
inherits an action of~\(G\) by naturality (which need not be strongly
continuous).

\begin{definition}
  \label{def:asymp_mor_G}
  Let \(A\) and~\(B\) be two \(G\)\nbd{}\Cstar{}algebras for a locally
  compact group~\(G\).  A \emph{\(G\)\nbd{}equivariant asymptotic
    morphism} from~\(A\) to~\(B\) is a \(G\)\nbd{}equivariant
  \Star{}homomorphism \(f\colon A\to\Asymp(B)\).  We write \(\llbracket
  A,B\rrbracket\) for the set of homotopy classes of
  \(G\)\nbd{}equivariant asymptotic morphisms from~\(A\) to~\(B\).  Here
  a homotopy is a \(G\)\nbd{}equivariant \Star{}homomorphism
  \(A\to\Asymp\bigl(\CONT([0,1],B)\bigr)\).
\end{definition}

The asymptotic algebra fits, by definition, into an extension
\[
\CONT_0(\R_+,B)\into \CONT_b(\R_+,B)\prto \Asymp(B).
\]
Notice that \(\CONT_0(\R_+,B) \cong \Cone(B)\) is contractible.  If
\(f\colon A\to\Asymp(B)\) is a \(G\)\nbd{}equivariant asymptotic
morphism, then we can use it to pull back this extension to an extension
\(\Cone(B) \into E\prto A\) in \(\CstarcatG\); the \(G\)\nbd{}action
on~\(E\) is automatically strongly continuous.  If~\(F\) is exact and
homotopy invariant, then \(F\bigl(\Sus^n(E)\bigr)\to
F\bigl(\Sus^n(A)\bigr)\) is an isomorphism for all \(n\ge1\) by
Proposition~\ref{pro:exact_longsequence}.  The evaluation map
\(\CONT_b(\R_+,B)\to B\) at some \(t\in\R_+\) pulls back to a morphism
\(E\to B\), and these morphisms for different~\(t\) are all homotopic.
Hence we get a well-defined map \(F\bigl(\Sus^n(A)\bigr) \cong
F\bigl(\Sus^n(E)\bigr)\to F\bigl(\Sus^n(B)\bigr)\) for each asymptotic
morphism \(A\to B\).  This explains how asymptotic morphisms are related
to exact homotopy functors.  This observation leads to the following
theorem:

\begin{theorem}
  \label{the:E_concrete}
  There are natural bijections
  \[
  \E^G_0(A,B) \cong \bigl\llbracket
  \Sus\bigl(A_\Comp \otimes \Comp(\ell^2\N)\bigr),
  \Sus\bigl(B_\Comp \otimes \Comp(\ell^2\N)\bigr)
  \bigr\rrbracket
  \]
  for all separable \(G\)\nbd{}\Cstar{}algebras \(A,B\).
\end{theorem}

An important step in the proof of Theorem~\ref{the:E_concrete} is the
\emph{\textsc{Connes}--\textsc{Higson} construction}, which to an
extension \(I\into E\prto Q\) in \(\Cstarsep\) associates an asymptotic
morphism \(\Sus(Q)\to I\).  A \(G\)\nbd{}equivariant generalisation of
this construction is discussed in
\cite{Thomsen:Asymptotic_Equivariant_KK}.  Thus any extension
in~\(\CstarsepG\) gives rise to an exact triangle \(\Sus(Q)\to I\to E\to
Q\) in~\(\Ecat^G\).

This also leads to the triangulated category structure of \(\Ecat^G\).
As for \(\KKcat^G\), we can define it using mapping cone triangles or
extension triangles \mdash{} both approaches yield the same class of
exact triangles.  The canonical functor \(\KKcat^G\to \Ecat^G\) is exact
because it evidently preserves mapping cone triangles.

Now that we have two bivariant homology theories with apparently very
similar formal properties, we must ask which one we should use.  It may
seem that the better exactness properties of \(\E\)\nbd{}theory raise it
above \(\KK\)-theory.  But actually, these strong exactness properties
have a drawback: for a general group~\(G\), the reduced crossed product
functor need not be exact, so that there is no guarantee that it extends
to a functor \(\Ecat^G\to\Ecat\).  Only full crossed products exist for
all groups by Proposition~\ref{pro:fullcross_exact}; the construction of
\(G\cross\blank\colon \Ecat^G\to\Ecat\) is the same as in
\(\KK\)-theory.  Similar problems occur with \(\minotimes\) but not with
\(\maxotimes\).

Furthermore, since \(\KKcat^G\) has a weaker universal property, it acts
on more functors, so that results about \(\KKcat^G\) have stronger
consequences.  A good example of a functor that is split-exact but
probably not exact is \emph{local cyclic cohomology} (see
\cites{Meyer:HLHA, Puschnigg:Diffeotopy}).  Therefore, the best practice
seems to prove results in \(\KKcat^G\) if possible.

Many applications can be done with either \(\E\) or \(\KK\), we hardly
notice any difference.  An explanation for this is the work of
\textsc{Houghton-Larsen} and \textsc{Thomsen}
(\cites{Houghton-Larsen-Thomsen, Thomsen:Asymptotic_Equivariant_KK}),
which describes \(\KK^G_0(A,B)\) in the framework of asymptotic
morphisms.  Recall that asymptotic morphisms \(A\to B\) generate
extensions \(\Cone(B)\into E\prto A\).  If this extension is
\(G\)\nbd{}equivariantly cp-split, then the projection map \(E\prto
A\) is a \(\KK^G\)-equivalence.  A \(G\)\nbd{}equivariant completely
positive contractive section for the extension exists if and only if
we can represent our asymptotic morphism by a continuous family of
\(G\)\nbd{}equivariant, completely positive contractions \(f_t\colon
A\to B\), \(t\in[0,\infty)\).

\begin{definition}
  \label{def:asymp_mor_cp}
  Let \(\llbracket A,B\rrbracket_{\mathrm{cp}}\) be the set of
  homotopy classes of asymptotic morphisms from~\(A\) to~\(B\) that
  can be lifted to a \(G\)\nbd{}equivariant, completely positive,
  contractive map \(A\to \CONT_b(\R_+,B)\); of course, we only use
  homotopies with the same kind of lifting.
\end{definition}

\begin{theorem}[\cite{Thomsen:Asymptotic_Equivariant_KK}]
  \label{the:KK_like_E}
  There are natural bijections
  \[
  \KK^G_0(A,B) \cong \bigl\llbracket
  \Sus\bigl(A_\Comp \otimes \Comp(\ell^2\N)\bigr),
  \Sus\bigl(B_\Comp \otimes \Comp(\ell^2\N)\bigr)
  \bigr\rrbracket_{\mathrm{cp}}
  \]
  for all separable \(G\)\nbd{}\Cstar{}algebras \(A,B\); the canonical
  functor \(\KKcat^G\to \Ecat^G\) corresponds to the obvious map that
  forgets the additional constraints.
\end{theorem}

\begin{corollary}
  \label{cor:KK_E_nuclear}
  If~\(A\) is a nuclear \Cstar{}algebra, then \(\KK_*(A,B)\cong
  \E_*(A,B)\).
\end{corollary}

\begin{proof}
  The \textsc{Effros}--\textsc{Choi} Lifting Theorem asserts that any
  extension of~\(A\) has a completely positive contractive section.
\end{proof}

In the equivariant case, the same argument yields \(\KK_*^G(A,B)\cong
\E_*^G(A,B)\) if~\(A\) is nuclear and~\(G\) acts properly on~\(A\) (see
also \cite{Skandalis:K-nuclear}).  It should be possible to weaken
properness to amenability here, but I am not aware of a reference for
this.

\section{The \textsc{Baum}--\textsc{Connes} assembly map for spaces and
  operator algebras}
\label{sec:BC_spaces_opalg}

The \textsc{Baum}--\textsc{Connes} conjecture is a guess for the
\(\K\)\nbd{}theory \(\K_*(\Cred G)\) of reduced group \Cstar{}algebras
\cite{Baum-Connes-Higson}.  We shall compare the approach of
\textsc{Davis} and \textsc{L\"uck} \cite{Davis-Lueck:Assembly} using
homotopy theory for \(G\)\nbd{}spaces and its counterpart in bivariant
\(\K\)\nbd{}theory formulated in~\cite{Meyer-Nest:BC}.  To avoid
technical difficulties, we \emph{assume that the group~\(G\) is
  discrete}.

The first step in the \textsc{Davis}--\textsc{L\"uck} approach is to
embed the groups of interest such as \(\K_*(\Cred G)\) in a
\emph{\(G\)\nbd{}homology theory}, that is, a homology theory on the
category of (spectra of) \(G\)\nbd{}CW-complexes.  For the
\textsc{Baum}--\textsc{Connes} assembly map we need a homology theory
for \(G\)\nbd{}CW-complexes with \(F_*(G/H) \cong \K_*(\Cred H)\).  This
amounts to finding a \(G\)\nbd{}equivariant spectrum with appropriate
homotopy groups \cite{Davis-Lueck:Assembly} and is the most difficult
part of the construction.  Other interesting invariants like the
algebraic \(\K\)- and \(\mathrm{L}\)\nbd{}theory of group rings can be
treated using other spectra instead.

In the world of \Cstar{}algebras, we cannot treat algebraic \(\K\)- and
\(\mathrm{L}\)\nbd{}theory; but we have much better tools to study
\(\K_*(\Cred G)\).  We do not need a \(G\)\nbd{}homology theory but a
homological functor on the triangulated category \(\KKcat^G\).  More
precisely, we need a homological functor that takes the value
\(\K_*(\Cred H)\) on \(\CONT_0(G/H)\) for all subgroups~\(H\).  The
functor \(A\mapsto \K_*(G\rcross A)\) works fine here because \(G\rcross
\CONT_0(G/H,A)\) is \textsc{Morita}--\textsc{Rieffel} equivalent to
\(H\rcross A\) for any \(H\)\nbd{}\Cstar{}algebra~\(A\).  The
corresponding assertion for full crossed products is known as
\emph{\textsc{Green}'s Imprimitivity Theorem}; reduced crossed products
can be handled similarly.  Thus a topological approach to the
\textsc{Baum}--\textsc{Connes} conjecture forces us to consider
\(\K_*(G\rcross A)\) for all \(G\)\nbd{}\Cstar{}algebras~\(A\), which
leads to the \emph{\textsc{Baum}--\textsc{Connes} conjecture with
  coefficients}.

\subsection{Assembly maps via homotopy theory}
\label{sec:proper_G_spaces}

Recall that a homology theory on pointed CW-complexes is determined by
its value on~\(\Sphere^0\).  Similarly, a \(G\)\nbd{}homology
theory~\(F\) is determined by its values \(F_*(G/H)\) on homogeneous
spaces for all subgroups \(H\subseteq G\).  This does not help much
because these groups \mdash{} which are \(\K_*(\Cred H)\) in the case of
interest \mdash{} are very hard to compute.

The idea behind assembly maps is to approximate a given homology theory
by a simpler one that only depends on \(F_*(G/H)\) for \(H\in\Fam\) for
some family of subgroups~\(\Fam\).  The \textsc{Baum}--\textsc{Connes}
assembly map uses the family of finite subgroups here; other families
like virtually cyclic subgroups appear in isomorphism conjectures for
other homology theories.  We now fix a family of subgroups~\(\Fam\),
which we assume to be closed under conjugation and subgroups.

A \emph{\(G,\Fam\)-CW-complex} is a \(G\)\nbd{}CW-complex in which the
stabilisers of cells belong to~\(\Fam\).  The \emph{universal
  \(G,\Fam\)-CW-complex} is a \(G,\Fam\)\nbd{}CW-complex \(\EG(G,\Fam)\)
with the property that, for any \(G,\Fam\)-CW-complex~\(X\) there is a
\(G\)\nbd{}map \(X\to \EG(G,\Fam)\), which is unique up to
\(G\)\nbd{}homotopy.  This universal property determines \(\EG(G,\Fam)\)
uniquely up to \(G\)\nbd{}homotopy.  It is easy to see that
\(\EG(G,\Fam)\) is \(H\)\nbd{}equivariantly contractible for any
\(H\in\Fam\).  Conversely, a \(G,\Fam\)-CW-complex with this property is
universal.

\begin{example}
  \label{exa:EG_Z}
  Let \(G=\Z\) and let~\(\Fam\) be the family consisting only of the
  trivial subgroup; this agrees with the family of finite subgroups
  because~\(G\) is torsion-free.  A \(G,\Fam\)-CW-complex is
  essentially the same as a CW-complex with a free cellular action
  of~\(\Z\).  It is easy to check that~\(\R\) with the action of~\(\Z\)
  by translation and the usual cell decomposition is a universal
  \(G,\Fam\)-CW-complex.
\end{example}

Given any \(G\)\nbd{}CW-complex~\(X\), the canonical map
\(\EG(G,\Fam)\times X\to X\) has the following properties:
\begin{itemize}
\item \(\EG(G,\Fam)\times X\) is a \(G,\Fam\)\nbd{}CW-complex;

\item if~\(Y\) is a \(G,\Fam\)\nbd{}CW-complex, then any \(G\)\nbd{}map
  \(Y\to X\) lifts uniquely up to \(G\)\nbd{}homotopy to a map \(Y\to
  \EG(G,\Fam)\times X\);

\item for any \(H\in\Fam\), the map \(\EG(G,\Fam)\times X\to X\) becomes
  a homotopy equivalence in the category of \(H\)\nbd{}spaces.

\end{itemize}
The first two properties make precise in what sense \(\EG(G,\Fam)\times
X\) is the best approximation to~\(X\) among \(G,\Fam\)-CW-complexes.

\begin{definition}
  \label{def:assembly_spaces}
  The \emph{assembly map} with respect to~\(\Fam\) is the map
  \(F_*\bigl(\EG(G,\Fam)\bigr) \to F_*(\pt)\) induced by the constant
  map \(\EG(G,\Fam) = \EG(G,\Fam)\times\pt\to\pt\).

  More generally, the assembly map with coefficients in a pointed
  \(G\)\nbd{}CW-complex (or spectrum)~\(X\) is the map
  \(F_*(\EG(G,\Fam)_+\wedge X)\to F_*(\Sphere^0\wedge X) = F_*(X)\)
  induced by the map \(\EG(G,\Fam)_+\to \pt_+=\Sphere^0\).
\end{definition}

In the stable homotopy category of pointed \(G\)\nbd{}CW-complexes (or
spectra), we get an exact triangle \(\EG(G,\Fam)_+\wedge X\to X\to N\to
\Sphere^1\wedge \EG(G,\Fam)_+\wedge X\), where~\(N\) is
\(H\)\nbd{}equivariantly contractible for each \(H\in\Fam\).  This means
that the domain of the assembly map \(F_*(\EG(G,\Fam)_+\wedge X)\) is
the \emph{localisation} of~\(F_*\) at the class of all objects that are
\(H\)\nbd{}equivariantly contractible for each \(H\in\Fam\).

Thus the assembly map is an isomorphism for all~\(X\) if and only if
\(F_*(N)=0\) whenever~\(N\) is \(H\)\nbd{}equivariantly contractible for
each \(H\in\Fam\).  Thus an isomorphism conjecture can be interpreted in
two equivalent ways.  First, it says that we can reconstruct the
homology theory from its restriction to \(G,\Fam\)\nbd{}CW-complexes.
Secondly, it says that the homology theory vanishes for spaces that are
\(H\)\nbd{}equivariantly contractible for \(H\in\Fam\).

\subsection{From spaces to operator algebras}
\label{sec:spaces_to_operator_algebras}

We can carry over the construction of assembly maps above to bivariant
\textsc{Kasparov} theory; we continue to assume~\(G\) discrete to
simplify some statements.  From now on, we let~\(\Fam\) be the family of
finite subgroups.  This is the family that appears in the
\textsc{Baum}--\textsc{Connes} assembly map.  Other families of
subgroups can also be treated, but some proofs have to be modified and
are not yet written down.

First we need an analogue of \(G,\Fam\)-CW-complexes.  These are
constructible out of simpler ``cells'' which we describe first, using
the \emph{induction functors}
\[
\Ind_H^G\colon \KKcat^H\to\KKcat^G
\]
for subgroups \(H\subseteq G\).  For a finite group~\(H\),
\(\Ind_H^G(A)\) is the \(H\)\nbd{}fixed point algebra of
\(\CONT_0(G,A)\), where~\(H\) acts by \(h\cdot f(g) = \alpha_h\bigl(
f(gh)\bigr)\).  For infinite~\(H\), we have
\begin{multline*}
  \Ind_H^G(A) = \{f\in \CONT_b(G,A)\mid
  \\\text{\(\alpha_h f(gh)=f(g)\) for all \(g\in G\), \(h\in H\), and
    \(gH\mapsto \norm{f(g)}\) is in \(\CONT_0(G/H)\)}\};
\end{multline*}
the group~\(G\) acts by translations on the left.

This construction is functorial for equivariant \Star{}homomorphisms.
Since it commutes with \Cstar{}stabilisations and maps split
extensions again to split extensions, it descends to a functor
\(\KKcat^H\to\KKcat^G\) by the universal property (compare
\textsection\ref{sec:construct_functors_universal}).

We also have the more trivial \emph{restriction functors}
\(\Res_G^H\colon \KKcat^G\to\KKcat^H\) for subgroups \(H\subseteq G\).
The induction and restriction functors are adjoint:
\[
\KK^G(\Ind_H^G A,B)\cong \KK^H(A,\Res_G^H B)
\]
for all \(A\inOb\KKcat^G\); this can be proved like the similar
adjointness statements in
\textsection\ref{sec:construct_functors_universal}, using the embedding
\(A\to \Res_G^H\Ind_H^G(A)\) as functions supported on \(H\subseteq G\)
and the correspondence \(\Ind_H^G\Res_G^H(A)\cong \CONT_0(G/H,A)\to
\Comp(\ell^2 G/H)\otimes A\sim A\).  It is important here that
\(H\subseteq G\) is an \emph{open} subgroup.  By the way, if
\(H\subseteq G\) is a \emph{cocompact} subgroup (which means finite
index in the discrete case), then \(\Res_G^H\) is the left-adjoint of
\(\Ind_H^G\) instead.

\begin{definition}
  \label{def:CI}
  We let \(\mathcal{CI}\) be the subcategory of all objects of
  \(\KKcat^G\) of the form \(\Ind_H^G(A)\) for \(A\inOb\KKcat^H\) and
  \(H\in\Fam\).  Let \(\genci\) be the smallest class in \(\KKcat^G\)
  that contains \(\mathcal{CI}\) and is closed under
  \(\KK^G\)-equivalence, countable direct sums, suspensions, and exact
  triangles.
\end{definition}

Equivalently, \(\genci\) is the \emph{localising subcategory} generated
by~\(\mathcal{CI}\).  This is our substitute for the category of
\((G,\Fam)\)-CW-complexes.

\begin{definition}
  \label{def:CC}
  Let \(\CC\) be the class of all objects of \(\KKcat^G\) with
  \(\Res_G^H(A)\cong 0\) for all \(H\in\Fam\).
\end{definition}

\begin{theorem}
  \label{the:CC_CI}
  If \(P\in\genci\), \(N\in\CC\), then \(\KK^G(P,N)=0\).  Furthermore,
  for any \(A\inOb\KKcat^G\) there is an exact triangle \(P\to A\to
  N\to\Sigma P\) with \(P\in\genci\), \(N\in\CC\).
\end{theorem}

Definitions \ref{def:CI}--\ref{def:CC} and Theorem~\ref{the:CC_CI} are
taken from~\cite{Meyer-Nest:BC}.  The map \(F_*(P)\to F_*(A)\) for a
functor \(F\colon \KKcat^G\to\Cat\) is analogous to the assembly map in
Definition~\ref{def:assembly_spaces} and deserves to be called the
\emph{\textsc{Baum}--\textsc{Connes} assembly map for~\(F\)}.

We can use the tensor product in \(\KKcat^G\) to simplify the proof of
Theorem~\ref{the:CC_CI}: once we have a triangle \(P_\C\to \C\to N_\C\to
\Sigma P_\C\) with \(P_\C\in\genci\), \(N_\C\in\CC\), then
\[
A\otimes P_\C \to A\otimes \C\to A\otimes N_\C \to
\Sigma A\otimes P_\C
\]
is an exact triangle with similar properties for~\(A\).  It makes no
difference whether we use \(\minotimes\) or \(\maxotimes\) here.  The
map \(P_\C\to\C\) in \(\KK^G(P_\C,\C)\) is analogous to the map
\(\EG(G,\Fam)\to\pt\).  It is also called a \emph{\textsc{Dirac}
  morphism} for~\(G\) because the \(\K\)\nbd{}homology classes of
\textsc{Dirac} operators on smooth spin manifolds provided the first
important examples \cite{Kasparov:Novikov}.

The two assembly map constructions with spaces and \Cstar{}algebras are
not just analogous but provide the same \textsc{Baum}--\textsc{Connes}
assembly map.  To see this, we must understand the passage from the
homotopy category of spaces to \(\KKcat\).  Usually, we map spaces to
operator algebras using the commutative \Cstar{}algebra \(\CONT_0(X)\).
But this construction is only functorial for \emph{proper} continuous
maps, and the functoriality is contravariant.  The assembly map for,
say, \(G=\Z\) is related to the non-proper map \(p\colon \R\to\pt\),
which does not induce a map \(\C\to \CONT_0(\R)\); even if it did, this
map would still go in the wrong direction.  The \emph{wrong-way
  functoriality} in \(\KK\) provides an element \(p_!\in
\KK_1(\CONT_0(\R),\C)\) instead, which is the desired \textsc{Dirac}
morphism up to a shift in the grading.  This construction only applies
to manifolds with a Spin\(^c\)-structure, but it can be generalised as
follows.

On the level of \textsc{Kasparov} theory, we can define another
functor from suitable spaces to \(\KKcat\) that is a \emph{covariant}
functor for \emph{all} continuous maps.  The definition uses a notion
of duality due to \textsc{Kasparov} \cite{Kasparov:Novikov} that is
studied more systematically in~\cite{Emerson-Meyer:Gysin}.  It
requires yet another version \(\RKK^G_*(X;A,B)\) of \textsc{Kasparov}
theory that is defined for a locally compact space~\(X\) and two
\(G\)\nbd{}\Cstar{}algebras \(A\) and~\(B\).  Roughly speaking, the
cycles for this theory are \(G\)\nbd{}equivariant families of cycles
for \(\KK_*(A,B)\) parametrised by~\(X\).  The groups
\(\RKK^G_*(X;A,B)\) are contravariantly functorial and homotopy
invariant in~\(X\) (for \(G\)\nbd{}equivariant continuous maps).

We have \(\RKK^G_*(\pt;A,B) = \KK^G_*(A,B)\) and, more generally,
\(\RKK^G_*(X;A,B) \cong \KK^G_*\bigl(A,\CONT(X,B)\bigr)\) if~\(X\) is
\emph{compact}.  The same statement holds for non-compact~\(X\), but the
algebra \(\CONT(X,B)\) is not a \Cstar{}algebra any more: it is an
inverse system of \Cstar{}algebras.

\begin{definition}[\cite{Emerson-Meyer:Gysin}]
  \label{def:Kasparov_dual}
  A \(G\)\nbd{}\Cstar{}algebra~\(P_X\) is called an \emph{abstract dual}
  for~\(X\) if, for all second countable locally compact
  \(G\)\nbd{}spaces~\(Y\) and all separable \(G\)\nbd{}\Cstar{}algebras
  \(A\) and~\(B\), there are natural isomorphisms
  \[
  \RKK^G(X\times Y;A,B) \cong \RKK^G(Y;P_X\otimes A,B)
  \]
  that are compatible with tensor products.
\end{definition}

Abstract duals exist for many spaces.  For trivial reasons, \(\C\) is an
abstract dual for the one-point space.  For a smooth manifold~\(X\) with
an isometric action of~\(G\), both \(\CONT_0(T^*X)\) and the algebra of
\(\CONT_0\)\nbd{}sections of the \textsc{Clifford} algebra bundle
on~\(X\) are abstract duals for~\(X\); if~\(X\) has a
\(G\)\nbd{}equivariant Spin\(^c\)-structure \mdash{} as in the example
of~\(\Z\) acting on~\(\R\) \mdash{} we may also use a suspension of
\(\CONT_0(X)\).  For a finite-dimensional simplicial complex with a
simplicial action of~\(G\), an abstract dual is constructed by
\textsc{Gennadi Kasparov} and \textsc{Georges Skandalis}
in~\cite{Kasparov-Skandalis:Buildings} and in more detail
in~\cite{Emerson-Meyer:Gysin}.  It seems likely that the construction
can be carried over to infinite-dimensional simplicial complexes as
well, but this has not yet been written down.

There are also spaces with no abstract dual.  A prominent example is the
\textsc{Cantor} set: it has no abstract dual, even for trivial~\(G\)
(see~\cite{Emerson-Meyer:Gysin}).

Let~\(\mathcal{D}\) be the class of all \(G\)\nbd{}spaces that admit a
dual.  Recall that \(X\mapsto \RKK^G(X\times Y;A,B)\) is a contravariant
homotopy functor for continuous \(G\)\nbd{}maps.  Passing to
corepresenting objects, we get a \emph{covariant} homotopy functor
\[
\mathcal{D}\to\KKcat^G,\qquad X\mapsto P_X.
\]
This functor is very useful to translate constructions from homotopy
theory to bivariant \(\K\)\nbd{}theory.  An instance of this is the
comparison of the \textsc{Baum}--\textsc{Connes} assembly maps in both
setups:

\begin{theorem}
  \label{the:BC_compare}
  Let~\(\Fam\) be the family of finite subgroups of a discrete
  group~\(G\), and let \(\EG(G,\Fam)\) be the universal
  \((G,\Fam)\)-CW-complex.  Then \(\EG(G,\Fam)\) has an abstract
  dual~\(P\), and the map \(\EG(G,\Fam)\to\pt\) induces a \textsc{Dirac}
  morphism in \(\KK^G_0(P,\C)\).
\end{theorem}

Theorem~\ref{the:BC_compare} should hold for all families of
subgroups~\(\Fam\), but only the above special case is treated in
\cites{Emerson-Meyer:Gysin,Meyer-Nest:BC}.

\subsection{The \textsc{Dirac}-dual-\textsc{Dirac} method and geometry}
\label{sec:Dirac-dual-Dirac}

Let us compare the approaches in \textsection\ref{sec:proper_G_spaces}
and \textsection\ref{sec:spaces_to_operator_algebras}!  The bad thing
about the \Cstar{}algebraic approach is that it applies to fewer
theories.  The good thing about it is that \textsc{Kasparov} theory is
so flexible that any canonical map between \(\K\)\nbd{}theory groups has
a fair chance to come from a morphism in \(\KKcat^G\) which we can
construct explicitly.

For some groups, the \textsc{Dirac} morphism in \(\KK^G(P,\C)\) is a
\(\KK\)\nbd{}equivalence:

\begin{theorem}[\textsc{Higson}--\textsc{Kasparov}
  \cite{Higson-Kasparov:Amenable}]
  \label{the:BC_amenable}
  Let the group~\(G\) be amenable or, more generally, a-T-menable.  Then
  the \textsc{Dirac} morphism for~\(G\) is a \(\KK^G\)-equivalence, so
  that~\(G\) satisfies the \textsc{Baum}--\textsc{Connes} conjecture
  with coefficients.
\end{theorem}

The class of groups for which the \textsc{Dirac} morphism has a
\emph{one-sided} inverse is even larger.  This is the point of the
\emph{\textsc{Dirac}-dual-\textsc{Dirac} method}.  The following
definition in~\cite{Meyer-Nest:BC} is based on a simplification of this
method:

\begin{definition}
  \label{def:dual_Dirac}
  A \emph{dual \textsc{Dirac} morphism} for~\(G\) is an element
  \(\eta\in\KK^G(\C,P)\) with \(\eta\circ D=\ID_P\).
\end{definition}

If such a dual \textsc{Dirac} morphism exists, then it provides a
section for the assembly map \(F_*(P\otimes A)\to F_*(A)\) for any
functor \(F\colon \KKcat^G\to\Cat\) and any \(A\inOb\KKcat^G\), so that
the assembly map is a split monomorphism.  Currently, we know no group
without a dual \textsc{Dirac} morphism.  It is shown
in~\cites{Emerson-Meyer1,Emerson-Meyer2,Emerson-Meyer3} that the
existence of a dual \textsc{Dirac} morphism is a \emph{geometric
  property} of~\(G\) because it is related to the invertibility of
another assembly map that only depends on the \emph{coarse} geometry
of~\(G\) (in the torsion-free case).

Instead of going into this construction, we briefly indicate another
point of view that also shows that the existence of a dual
\textsc{Dirac} morphism is a geometric issue.  Let~\(P\) be an abstract
dual for some space~\(X\) (like \(\EG(G,\Fam)\)).  The duality
isomorphisms in Definition~\ref{def:Kasparov_dual} are determined by two
pieces of data: a \emph{\textsc{Dirac} morphism} \(D\in\KK^G(P,\C)\) and
a \emph{local dual \textsc{Dirac} morphism} \(\Theta\in\RKK^G(X;\C,P)\).
The notation is motivated by the special case of a
Spin\(^c\)-manifold~\(X\) with \(P=\CONT_0(X)\), where~\(D\) is the
\(\K\)\nbd{}homology class defined by the \textsc{Dirac} operator
and~\(\eta\) is defined by a local construction involving pointwise
\textsc{Clifford} multiplications.  If \(X=\EG(G,\Fam)\), then it turns
out that \(\eta\in\KK^G(\C,P)\) is a dual \textsc{Dirac} morphism if and
only if the canonical map \(\KK^G(\C,P)\to\RKK^G(X;\C,P)\) maps
\(\eta\mapsto\Theta\).  Thus the issue is to \emph{globalise} the local
construction of~\(\Theta\).  This is possible if we know, say,
that~\(X\) has non-positive curvature.  This is essentially how
\textsc{Kasparov} proves the \textsc{Novikov} conjecture for fundamental
groups of non-positively curved smooth manifolds in
\cite{Kasparov:Novikov}.

\end{document}